\documentclass[12pt,leqno]{article}
\usepackage{amsfonts}
\usepackage{amsmath, amsthm, amsfonts, amssymb, color}
\usepackage{mathrsfs}
\usepackage{color}
\usepackage[notcite,notref]{showkeys}

\theoremstyle{definition}

\newcommand{\scr}[1]{\mathscr #1}
\definecolor{wco}{rgb}{0.5,0.2,0.3}

\numberwithin{equation}{section} \theoremstyle{remark}

\newcommand{\ua}{\uparrow}

\title{{\bf  Image Dependent Conditional McKean-Vlasov SDEs  for Measure-Valued Diffusion Processes}\footnote{Supported in
 part by  NNSFC (11771326, 11831014) and the DFG through the CRC 1283.} }
\author{
{\bf    Feng-Yu Wang$^{a),b)}$  }\\
\footnotesize{$^{a)}$ Center for Applied Mathematics, Tianjin University, Tianjin 300072, China }\\
 \footnotesize{ $^{b)}$ Department of Mathematics,
Swansea University, Singleton Park, SA2 8PP, United Kingdom}\\
\footnotesize{  wangfy@tju.edu.cn, F.-Y.Wang@swansea.ac.uk}}
\begin{document}
\allowdisplaybreaks
\def\R{\mathbb R}  \def\ff{\frac} \def\ss{\sqrt} \def\B{\mathbf
B}
\def\N{\mathbb N} \def\kk{\kappa} \def\m{{\bf m}}
\def\ee{\varepsilon}\def\ddd{D^*}
\def\dd{\delta} \def\DD{\Delta} \def\vv{\varepsilon} \def\rr{\rho}
\def\<{\langle} \def\>{\rangle} \def\GG{\Gamma} \def\gg{\gamma}
  \def\nn{\nabla} \def\pp{\partial} \def\E{\mathbb E}
\def\d{\text{\rm{d}}} \def\bb{\beta} \def\aa{\alpha} \def\D{\scr D}
  \def\si{\sigma} \def\ess{\text{\rm{ess}}}
\def\beg{\begin} \def\beq{\begin{equation}}  \def\F{\scr F}
\def\Ric{\mathcal Ric} \def\Hess{\text{\rm{Hess}}}
\def\e{\text{\rm{e}}} \def\ua{\underline a} \def\OO{\Omega}  \def\oo{\omega}
 \def\tt{\tilde}
\def\cut{\text{\rm{cut}}} \def\P{\mathbb P} \def\ifn{I_n(f^{\bigotimes n})}
\def\C{\scr C}      \def\aaa{\mathbf{r}}     \def\r{r}
\def\gap{\text{\rm{gap}}} \def\prr{\pi_{{\bf m},\varrho}}  \def\r{\mathbf r}
\def\Z{\mathbb Z} \def\vrr{\varrho} \def\ll{\lambda}
\def\L{\scr L}\def\Tt{\tt} \def\TT{\tt}\def\II{\mathbb I}
\def\i{{\rm in}}\def\Sect{{\rm Sect}}  \def\H{\mathbb H}
\def\M{\scr M}\def\Q{\mathbb Q} \def\texto{\text{o}} \def\LL{\Lambda}
\def\Rank{{\rm Rank}} \def\B{\scr B} \def\i{{\rm i}} \def\HR{\hat{\R}^d}
\def\to{\rightarrow}\def\l{\ell}\def\iint{\int}
\def\EE{\scr E}\def\Cut{{\rm Cut}}\def\W{\mathbb W}
\def\A{\scr A} \def\Lip{{\rm Lip}}\def\S{\mathbb S}
\def\BB{\scr B}\def\Ent{{\rm Ent}} \def\i{{\rm i}}\def\itparallel{{\it\parallel}}
\def\g{{\mathbf g}}\def\Sect{{\mathcal Sec}}\def\T{\mathcal T}\def\V{{\bf V}}
\def\PP{{\bf P}}\def\HL{{\bf L}}\def\Id{{\rm Id}}\def\f{{\bf f}}\def\cut{{\rm cut}}
\def\L{\scr L}

\def\BL{\scr A}

\maketitle

\begin{abstract} We consider a special class of mean field SDEs with common noise which depend on the image of  the solution (i.e. the conditional distribution given noise).  The strong well-posedness is derived under  a monotone condition which is weaker than those used in the literature of mean field games, the Feynman-Kac formula is established to solve Schr\"ordinegr type PDEs on $\scr P_2$, and the ergodicity is proved for a class of measure-valued diffusion processes. 
\end{abstract} \noindent
 \vskip 2cm

\section{Introduction}

Let $\scr P_2$ be the space of all probability measures $\mu$ on $\R^d$ such that
 $$ \|\mu\|_2:=\bigg(\int_{\R^d} |x|^2\mu(\d x)\bigg)^{\ff 1 2 } <\infty,$$ where $|\cdot|$ is the norm in $\R^d$. We will use $\|\cdot\|$ to denote the operator norm of a matrix or
  linear operator, and use $\|\cdot\|_{HS}$ to stand for the Hilbert-Schmidt norm.
   It is well known that  $\scr P_2$ is a Polish space under the Wasserstein distance
 $$\W_2(\mu,\nu):= \inf_{\pi\in \C(\mu,\nu)} \bigg(\int_{\R^d\times\R^d} |x-y|^2\pi(\d x,\d y)\bigg)^{\ff 1 2},$$
 where $\C(\mu,\nu)$ is the set of couplings for $\mu$ and $\nu$.

 Since 1996 when Albeverio, Kondratiev and R\"ockner \cite{AKR} introduced the intrinsic derivative  on the configuration space over manifolds,  diffusion processes  on the space of discrete Radon measures have been investigated by using Dirichlet forms,
 see \cite{KLV} and references within.
 This  derivative   provides a natural Riemannian structure on the Wasserstein space $(\scr P_2,\W_2)$, see Subsection 1.2 below.

To develop stochastic analysis and applications on this space, we intend to construct diffusion processes generated by second order differentiable operators and solve the associated PDEs on $\scr P_2$. Below we first recall the intrinsic/Lions derivative on $\scr P_2$.

According to \cite{AKR},  let   $L^2(\R^d\to\R^d;\mu)$ be the tangent space of $\scr P_2$ at point $\mu\in \scr P_2$, and define the directional derivative by
 $$D_\phi f(\mu):=  \lim_{\vv\downarrow 0} \ff{f(\mu\circ ({\rm Id}+\vv\phi)^{-1})-f(\mu)}{\vv},\ \ \phi\in L^2(\R^d\to\R^d;\mu).$$
 When $\phi\mapsto D_\phi f(\mu)$ is a bounded linear functional on  $L^2(\R^d\to\R^d;\mu)$, or equivalently the map
 \beq\label{*D0} L^2(\R^d\to \R^d;\mu)\ni \phi \mapsto f(\mu\circ({\rm Id} +\phi)^{-1})\end{equation}
is Gateaux differentiable at $\phi=0$,
  there exists a unique element $Df(\mu)\in L^2(\R^d\to\R^d;\mu)$ such that
$$\<D f(\mu), \phi\>_{L^2(\mu)}= D_\phi f(\mu),\ \ \phi\in L^2(\R^d\to\R^d;\mu).$$ In this case,  we call $f$ intrinsically differentiable at $\mu$ with derivative $Df(\mu)$.  According to Lions  (see \cite{Card}),  if $Df(\mu)$ exists and
\beq\label{*D1}\lim_{\mu(|\phi|^2)\to 0} \ff{f(\mu\circ ({\rm Id}+\phi)^{-1})-f(\mu)- D_\phi f(\mu)}{\ss{\mu(|\phi|^2)}}= 0,\end{equation} i.e. the map in \eqref{*D0} is Fr\'echet differentiable at  $\phi=0$,  we call $f$ $L$-differentiable at $\mu\in \scr P_2$.
If $f$ is $L$-differentiable at any $\mu\in \scr P_2$, we call it $L$-differentiable. Note that $Df(\mu)$ is a $\mu$-a.e. defined  $\R^d$-valued function. Let $\{Df(\mu)\}_i$ be its $i$-th component for $1\le i\le d.$

In this paper, we investigate diffusion processes and applications on the Wasserstein space $\scr P_2$.   Let $m\ge 1$, and let
$$b: [0,\infty)\times \R^d\times\scr P_2\to \R^d,\ \ \si: [0,\infty)\times \R^d\times\scr P_2\to \R^{d}\otimes\R^m$$ be measurable such that
$|b(t,\cdot,\mu)|+\|\si(t,\cdot,\mu)\|_{HS}^2\in L^1(\mu)$ for any $(t,\mu)\in [0,\infty)\times \scr P_2$.
We consider the following time-dependent second order differential operators on $\scr P_2$:
\beq\label{HA1} \beg{split}   \BL_t  f (\mu):= &\ \ff 1 2\int_{\R^d\times \R^d}  \big\<\si (t,y,\mu)\si(t,z,\mu)^*, D^2 f(\mu)(y,z)\big\> \mu(\d y)\mu(\d z) \\
& + \int_{\R^d} \Big(\ff 1 2   \big\< (\si\si^*)(t,y,\mu),  \nn\{Df(\mu)\}(y)\big\> +\big\<b(t,y,\mu), Df(\mu)(y)\big\> \Big)\mu(\d y),
\end{split}\end{equation}
where   $\<\cdot,\cdot\>$ is the inner product on $\R^d$ or $\R^d\otimes\R^d.$
We  also consider the    following  extension   of $\BL_t$ on $\R^d\times \scr P_2$:
\beq\label{HA2} \beg{split}   \tt \BL_t f(x,\mu):= &\ \BL_t f(x,\cdot)(\mu) +    \ff 1 2  \big\<\si (t,x,\mu)\si(t,x,\mu)^*, \nn^2 f(x,\mu) \big\> +\big\<b(t,x,\mu), \nn f(x,\mu)\big\>\\
& +  \int_{\R^d}   \big\<(D\nn f)(x,\mu) (y), \si(t,y,\mu)\si(t,x,\mu)^*\big\>  \mu(\d y).
\end{split}\end{equation}

To present   reasonable pre-domains of  $\BL_t$ and $\tt\BL_t$,   we introduce below some classes    of $L$-differentiable functions.
\beg{enumerate} \item[(1)]
We write $f\in C^{1}(\scr P_2)$, if $f$ is $L$-differentiable and the derivative has a $\mu$-version $Df(\mu)(x)$  which is jointly continuous in $(\mu,x)\in \scr P_2\times\R^d$. If moreover $Df(\mu)(x)$ is bounded in
$(x,\mu)\in \R^d\times\scr P_2$, we denote $f\in C^{1}_b(\scr P_2).$
\item[(2)] We write $f\in C^{(1,1)}(\scr P_2)$, if $f\in C^{1}_b(\scr P_2)$ and  $Df(\mu)(x)$ is differentiable in $x$ such that  the $\R^d\otimes\R^d$-valued function
$$\nn \{Df(\mu)\}(x) := \big(\pp_{x_j} \{Df(\mu)(x)\}_i\big)_{1\le i,j\le d}$$ is   jointly continuous in $(\mu,x)\in \scr P_2\times\R^d$. If moreover $Df(\mu)(x)$ and $\nn\{Df(\mu)\}(x)$ are bounded in
$(x,\mu)\in \R^d\times\scr P_2$,  we denote $f\in C_b^{(1,1)}(\scr P_2).$
\item[(3)] We write
 $f\in C^{2}(\scr P_2)$, if $f\in C^{(1,1)}(\scr P_2)$ and $Df(\mu)(x)$ is $L$-differentiable in $\mu$ such that the $\R^d\otimes\R^d$-valued function
$$ D^2f(\mu)(x,y) := \big(\big\{D[Df(\mu)(x)]_i (y)\big\}_j\big)_{1\le i,j\le d}$$   is  jointly continuous in $(\mu,x,y)\in \scr P_2\times\R^d\times\R^d.$ If moreover $f\in C_b^{(1,1)}(\scr P_2)$ and  $D^2f(\mu)(x,y)$  is bounded in
$(x,y,\mu)\in \R^d\times\R^d\times\scr P_2$, we denote $f\in C^{2}_b(\scr P_2)$.
 \item[(4)] We write $f\in C^{2,2}(\R^k\times\scr P_2)$ for some $k\ge 1$,  if $f$ is a continuous function on $\R^k\times \scr P_2$ such that
$f(\cdot,\mu)\in C^2(\R^k)$ for $\mu\in\scr P_2$, $f(x,\cdot)\in C^{2}(\scr P_2)$ for $x\in \R^k$,
$$(D\nn f)(x,\mu)(y):= \big(\big\{D[\pp_{x_i} f(x,\mu)]\big\}_j\big)_{1\le i,\j\le d}\in \R^d\otimes\R^d$$ exists, and the derivatives
$$\nn f(x,\mu), \nn^2f(x,\mu), Df(x,\mu)(y), (D\nn f)(x,\mu)(y), \nn\{Df(x,\mu)(\cdot)\}(y), D^2f(x,\mu)(y,z)$$ are bounded and jointly
continuous in the corresponding arguments.
\end{enumerate}

\paragraph{Example 1.1.} For any $p\ge 1$, consider the following class of cylindrical functions
\beq\label{FC} \beg{split}\F C_b^p(\scr P_2):= \big\{&f(\mu):=g(\mu(h_1),\cdots, \mu(h_n)):\\
&\quad  n\ge 1,
   g\in C_b^p(\R^n),  h_i\in C_b^p(\R^d), 1\le i\le n\big\}.\end{split}\end{equation}  When $p=2$,  such a function is in the class $C_b^{2}(\scr P_2)$ with
   \beq\label{LDD} \beg{split} & D f(\mu)(x)= \sum_{i=1}^n (\pp_ig )(\mu(h_1),\cdots, \mu(h_n))\nn h_i(x),\\
 &  D^2 f(\mu)(x,y)= \sum_{i,j=1}^n (\pp_i\pp_jg )(\mu(h_1),\cdots, \mu(h_n))\{\nn h_i(x)\}\otimes\{\nn h_j(y)\},\end{split}\end{equation}
 where $\{\nn h_i(x)\}\otimes\{\nn h_j(y)\}\in \R^d\otimes\R^d $ is defined as
 $$\big(\{\nn h_i(x)\}\otimes\{\nn h_j(y)\}\big)_{kl}= \{\pp_k h_i(x) \}\pp_l h_j(y),\ \ 1\le k,l\le d, x, y\in\R^d.$$
Moreover, $f\in C^{2,2}(\R^d\times \scr P_2)$ if $f(x,\mu)= g(x,\mu(h_1),\cdots, \mu(h_n))$ for some $n\ge 1$, $g\in C_b^2(\R^{n+d})$ and $\{h_i\}_{1\le i\le n}\subset C_b^2(\R^d)$.

\

We will construct the $\BL_t$-diffusion process  by solving    the following SDE  on $\R^d$:
\beq\label{E1} \d X_{s,t}^{x,\mu}= b(t,X_{s,t}^{x,\mu},  \LL_{s,t}^\mu)\d t +\si(t,X_{s,t}^{x,\mu},\LL_{s,t}^\mu)\d W_t,\ \LL_{s,t}^\mu:= \mu\circ(X_{s,t}^{\cdot,\mu})^{-1},\ t\ge s, X_{s,s}^{x,\mu}=x, \end{equation}
where $W_t$ is the $m$-dimensional Brownian motion on a complete filtration probability space $(\OO,\{\F_t\}_{t\ge 0},\P)$, $(s,x,\mu)\in [0,\infty)\times \R^d\times \scr P_2.$  Since this SDE depends on the image of   solutions, we call it   image dependent SDE.

It turns out that  the solution of \eqref{E1} for $s=0$ gives rise to a strong solution to the following conditional McKean-Vlasov  SDE
arising from mean field games:
\beq\label{E''}  \d X_{t}= b(t,X_{t}, \L_{X_{t}| W})\d t +\si(t,X_{t},\L_{X_{t}| W})\d W_t,\    \ \L_{X_{0}}= \mu\in \scr P_2, \end{equation}
where $\L_{\xi}$ and $\L_{\xi|W}$  denote the distribution and the conditional distribution given $\{W_t:t\ge 0\}$ for a random variable $\xi$. More precisely, $X_t=X_{0,t}^{X_0,\mu}$, see the proof of Corollary \ref{C2.2} below.  So, the measure-valued process $\LL_{0,t}$ in \eqref{E1}  is indeed the conditional distribution of $X_t$ given the noise $W$, and we may  call \eqref{E1} an image dependent conditional McKean-Vlasov SDE as in the title of the paper. 

The weak solution to \eqref{E''} has been investigated by using mean filed games with common noise.  More precisely, let $\{x_i\}_{i\ge }$ be a sequence of points in $\R^d$ such that 
$$\lim_{n\to\infty} \ff 1 n \sum_{i=1}^n \dd_{x_i} =\mu\  \text{ weakly},$$
consider the SDEs
$$\d X_t^{n,i}= b\Big(t,X_t^{n,i},\ff 1 n \sum_{j=1}^n \dd_{X_{t}^{n,j}}\Big)\d t +\si\Big(t,X_t^{n,i},\ff 1 n \sum_{j=1}^n \dd_{X_{t}^{n,j}}\Big)\d W_t,\ \ X_0^{n,i}=x_i, 1\le i\le n.$$
Then under reasonable conditions, when $n\to\infty$ the law of  $(X_t^{n,1})_{t\ge 0}$ converges weakly to a probability measure on the path space 
$C([0,\infty);\R^d)$ which solves \eqref{E''} weakly.  See \cite{CD18, CDL16, CDLL19,  DLR19, HSS19} for the study of a more general model than \eqref{E''} where an additional independent Brownian noise $W_t^0$ is included:
\beq\label{E'1}  \d X_{t}= b(t,X_{t}, \L_{X_{t}| W})\d t +\si(t,X_{t},\L_{X_{t}| W})\d W_t+ \si^0(t, X_t, \L_{X_t|W})\d W_t^0\end{equation} for $ \L_{X_{0}}= \mu\in \scr P_2,$
where $\si^0$ takes values in $\R^d\otimes\R^l$ if $W_t^0$ is $l$-dimensional. The study of this type SDE using mean field games goes back to the poineering works of Lasry and Lions \cite{LL06a, LL06b, LL07} and Huang, Malham\'e and Caines \cite{HMC06a, HMC06b}, see the nice monograph \cite{CD18} for a   theory of mean field games with common noises and applications.
 In this paper, we will study the strong solutions and applications  of \eqref{E1} (hence, \eqref{E''}) in a straightforward way under reasonably weaker conditions on the coefficients. 
 
In the remainder of this section, we first summarize  the main results of the paper, then present a link of the present model to the Brownian motion on   $\scr P_2$ for further study, and finally introduce some previous work for analysis on the Wasserstein space.

\subsection{Summary of main results}
 \paragraph{Existence and uniqueness.} Under a monotone condition, Theorem \ref{T2.1} ensures  the existence, uniqueness and moment estimates of solutions to the image SDE \eqref{E1}, and  that the unique solution is   the diffusion processes generated by $\BL_t$ on $\scr P_2$ and $\tt\BL_t$ on $\R^d\times \scr P_2$ respectively. As a consequence, the strong well-posedness is derived for the conditional distribution dependent SDE \eqref{E''}.  Our monotone condition is weaker than those  in \cite{CD18, CDL16} but incomparable with those of \cite{HSS19}, see Remark 2.1 below for details. 
  
\paragraph{Feynman-Kac formula.} By using the  diffusion process
$(X_{s,t}^{x,\mu},\LL_{s,t}^\mu)$,  Theorem \ref{T4.1} solves  the following PDE for $U$ on $[0,T]\times \R^d\times \scr P_2$ :   \beq\label{PDE}  \beg{split} & \pp_t U(t,x,\mu) + \tt\BL_t U(t,x,\cdot)(\mu) + (VU)(t,x,\mu) +  F(t,x,\mu)=0,\\
&\qquad  \  \ U(T,x,\mu)=\Phi(x,\mu), \ \ (t,x,\mu)\in [0,T]\times  \R^d\times \scr P_2,\end{split}\end{equation} where $T>0$ is a fixed time, $\Phi$ is a function on $\R^d\times \scr P_2$, and $V,F$ are functions on   $[0,T]\times\R^d\times  \scr P_2$. When $\Phi, F$ and $V$ do not depend on $x\in\R^d$, this PDE reduces to \beq\label{PDE2}  \beg{split} & \pp_t U(t,\mu) + \BL_t U(t,\cdot)(\mu) + (VU)(t,\mu) +  F(t,\mu)=0,\\
&\qquad  \  \ U(T,\mu)=\Phi(\mu), \ \ (t,\mu)\in [0,T]\times  \scr P_2.\end{split}\end{equation} When $V=0$ these two SPDEs are included as a special case by the Master equations studied in the literature of mean field games under stronger assumptions on $b$ and $\si$, see  Remark 3.1 below for details. 

\paragraph{Exponential ergodicity and structure of invariant probability measures.}Let $b$ and $\si$ do not depend on $t$. Under a dissipativity condition,
Theorem \ref{T3.1} provides  the exponential convergence rate of the   diffusion process $(X_t^{x,\mu},\LL_t^\mu):=(X_{0,t}^{x,\mu},\LL_{0,t}^\mu)$    to its unique invariant probability measure $\tt\Pi$. Consequently, the diffusion process $\LL_t^\mu$ converges at the same rate to  the invariant probability measure $\Pi:=\tt\Pi(\R^d\times\cdot)$.

Moreover,  let $b_0(x)= b(x,\dd_x), \si_0(x)= \si(x,\dd_x)$, and let $\mu_0$ be the unique invariant probability measure for the classical SDE
\beq\label{E0} \d X_t = b_0(X_t)\d t+ \si_0(X_t)\d W_t.\end{equation}
By Theorem \ref{T4.2}, $\tt\Pi$ and $\Pi$ have the representations
\beq\label{E01} \tt\Pi (\d x,\d \mu)= \mu_0(\d x) \dd_{\dd_x}(\d\mu),\ \ \Pi= \int_{\R^d} \dd_{\dd_x}  \mu_0(\d x),\end{equation}  where $\dd_{\dd_x}$ is the Dirac measure at point $\dd_x\in \scr P_2$.
This structure describes an asymptotic collision property of the diffusion process $\LL_t^\mu$: starting from   any probability measure $\mu\in \scr P_2$, the measure-valued process  eventually decays to a Dirac random variable, for which the whole mass focus on a single random point.

 \subsection{Some related studies}
\paragraph{Brownian motion on $\scr P_2$.} A Riemannian structure has been introduced in \cite{AS} on the Wasserstein space $(\scr P_2,\W_2)$.
With the intrinsic/Lions derivative, this space is  an infinite-dimensional Riemannian manifold with gradient $D$ and  Riemannian metric $\<\cdot,\cdot\>_{L^2(\mu)}$ on the tangent space $L^2(\R^d\to\R^d;\mu)$; that is,
$\W_2$ is the Riemannian distance induced by $D$.

As in the  finite-dimensional  Riemannian setting, we  introduce the square field
$$\GG(f,g)(\mu) :=\int_{\R^d} \<Df(\mu)(x), Dg(\mu)(x)\> \mu(\d x),\ \ f,g\in C_b^1(\scr P_2), $$
and the Laplace operator
$$\DD f(\mu):=  \int_{\R^d} {\rm tr}\big\{D^2f(\mu)(x,x)\big\}\mu(\d x),\ \ f\in C^2(\scr P_2).$$ Then by the chain rule we have
 $$\GG(f,g)= \ff 1 2\big\{ \DD (fg)- f\DD g- g\DD f\big\},\ \ f,g\in  C^2(\scr P_2).$$
This structure can be  easily extended to the Wasserstein space $\scr P_2(M)$  over a Riemannian manifold $M$. Note that when $M$ is compact we have $ \scr P_2(M)=\scr P(M)$, the space of all probability measures on $M$.

To develop stochastic analysis on $\scr P_2$, it is  interesting to  construct the Brownian motion, i.e. the diffusion process generated by $\ff 1 2\DD$;
or more generally, to construct diffusion processes on $\scr P_2$ with square field $\GG$. This  is the main motivation of \cite{RS} introduced in the next subsection.

Below we explain that when $\si\si^*={\rm Id}$ and $\mu=\dd_x$ is a Dirac measure at some point $x\in\R^d$, the process $(\LL_{0,t}^\mu)_{t\ge 0}$ is such a diffusion process.
Indeed, it is easy to check that the square field of the $\scr A_t$-diffusion process
is
\beg{align*} &\GG_t(f,g)(\mu):= \big\{\scr A_t (fg)(\mu) - f \scr A_t g -g \scr A_t \big\}(\mu)\\
&=  \int_{\R^d\times \R^d} \big\<\si(t,x,\mu)^* D f(\mu)(x),  \si(t,y,\mu)^* D f(\mu)(y)\big\>\mu(\d x)  \mu(\d y)\ \ f,g\in C_b^2(\scr P_2),\mu\in\scr P_2.\end{align*}
In particular, when $\si\si^*={\rm Id}$,  we have
$$\GG_t(f,g)(\mu)= \GG(f,g)(\mu),\ \  \mu\in \scr P_{2}^0:= \{\dd_x: x\in\R^d\}.$$ Since when $\mu=\dd_x$ for some $x\in \R^d$,     $\LL_{s,t}^{\mu}=\dd_{X_{s,t}^{x,\dd_x}}$ is a diffusion process on $\scr P_{2}^0$,
  Theorem \ref{T2.1}(2) below implies that $(\LL_{s,t}^\mu)_{t\ge s}$ for $\mu\in \scr P_2^0$   is a diffusion process with square field $\GG$.  However, this does not hold for $\mu\notin \scr P^0_{2}.$

\paragraph{Measure-valued diffusion processes.} Measure-valued diffusion processes  have been  constructed using Dirichlet forms.  Let $\scr P(\mathbb S^1)$ be the space of all  probability measures  on the unit circle $\mathbb S^1$.
A family of probability measures $\{\P_\bb\}_{\bb>0}$ on $\scr P(\mathbb S^1)$, called $``$entropic measures" with inverse temperature $\bb>0$,    have  been constructed by von Renesse and Sturm \cite{RS}  such that for each $\bb>0$, the bilinear form
$$\EE(f,g):=\int_{\scr P_2(\mathbb S^1)} \<Df(\mu), Dg(\mu)\>_{L^2(\mu)} \P_\bb(\d\mu)$$   gives   a symmetric Dirichlet form on $L^2(\Pi_\bb)$, which  refers to a
  $\P_\bb$-a.e. starting diffusion process    on $\scr P(\mathbb S^1)$.    See also \cite{Shao} for a different   Dirichlet form on $\scr P([0,1])$ with square field $\GG$. The construction of Dirichlet forms  in these papers  heavily relies on the
one-dimensional property.
See also \cite{KLV,RW19,W19,WZ19} and references within for the study of different type measure-valued diffusion processes using Dirichlet forms.

Next, following the idea of Konarovskyi (see e.g.\cite{ [12]}), \cite{KvR19,M18} constructed another type of  diffusion process on $\scr P([0,1])$  by  taking the limit as $N \to\infty$  
of a   system with $N$ coalescing and mass-carrying particles.  The generator $\L_t$  of the process has the formulation \cite[Theorem A.3]{M18}
$$\scr L_t f(\mu)= \ff 1 2 \E\int_0^1 \Big[\ff{\{Df(\mu)\}'(\xi_t(u))}{m_t(u)}+ D^2f(\mu)(\xi_t(u),\xi_t(u)) \Big]\d u,\ \ f\in C_b^2(\scr P([0,1])), t\ge 0,$$
where $\{(\xi_t(u))_{t\ge 0}:  u\in [0,1]\}$ is a family of continuous martingales with $ \xi\in D((0,1); C([0,\infty)) $ and quadratic variations \cite[Proposition 5.7]{M18}
$$ \<\xi(u),\xi(u)\>_t=\int_0^t\ff{\d s}{m_s(u)},\ \ m_s(u)=\int_0^1 1_{\{\xi_s(u)\ne \xi_s(v)\}} \d v.$$ Recently, the study of this type Wasserstein diffusion processes have been used in \cite{M20}  to solve a class of  Fokker-Planck equations on the interval  driven by  an infinite-dimensional noise.

\section{Image dependent SDE  and   diffusion processes on $\scr P_2$}

We will  construct the $\BL_t$-diffusion process by solving the image SDE \eqref{E1}.   In general,
we allow the coefficients
$$b: \OO\times [0,\infty)\times\scr P_2\to\R^d,\ \  \si: \OO\times [0,\infty)\times\scr P_2\to\R^d\otimes\R^m$$  to be random but progressively  measurable with respect to the filtration $\F_t$.
We first present the definition of solution.
\beg{defn} Let  $(s,\mu)\in [0,\infty)\times \scr P_2$. A family of adapted processes $\{(X_{s,t}^{x,\mu})_{t\ge s}:  x\in \R^d\}$ is called a solution to \eqref{E1}, if the following conditions     hold $\P$-a.s.:
   \beg{enumerate} \item[$(a)$]    $X_{s,t}^{x,\mu}$ is   continuous in $t\in [s,\infty)$ and measurable in $x\in\R^d$;
 \item[$(b)$]   $\LL_{s,t}^\mu:=\mu\circ (X_{s,t}^{\cdot,\mu})^{-1}\in \scr P_2$      is continuous in $t\ge s$;
 \item[$(c)$]    $ \E   \int_s^t  (|b(r, X_{s,r}^{x,\mu},\LL_{s,r}^\mu)|+ \|\si (r, X_{s,r}^{x,\mu}, \LL_{s,r}^\mu) \|_{HS}^2)\d r <\infty$ and
 $$X_{s,t}^{x,\mu}= x +\int_s^t b(r, X_{s,r}^{x,\mu}, \LL_{s,r}^\mu) \d r + \int_s^t \si (r, X_{s,r}^{x,\mu}, \LL_{s,r}^\mu) \d W_r,\ \ t\ge s, x\in\R^d.$$ \end{enumerate}  The image SDE  \eqref{E1} is called well-posed, if  it has a unique solution for any $(s,\mu)\in [0,\infty)\times\scr P_2$.
\end{defn}

To ensure the well-posedness of \eqref{E1}, we make the following assumption on $b$ and $\si$.

\beg{enumerate}\item[{\bf (A)}]   The progressively measurable coefficients $b(t,x,\mu)$ and $\si(t,x,\mu)$ are continuous in $ (x,\mu)\in \R^d\times\scr P_2$,   there exists
   $K\in L_{loc}^q([0,\infty)\to [0,\infty))$ for some $q>1$ such that $\P$-a.s. for any $t\ge 0,$
\beq\label{A1} |b(t,x,\mu)|^2+ \|\si(t,x,\mu)\|_{HS}^2 \le K(t) \big(1+|x|^2 + \|\mu\|_2^2\big),\ \ (x,\mu)\in \R^d\times\scr P_2,\end{equation}
\beq\label{A2} \beg{split}& 2\<b(t,x,\mu)- b(t,y,\nu), x-y\>^++ \|\si(t,x,\mu)-\si(t,y,\nu)\|_{HS}^2\\
&\qquad \le K(t) \big(|x-y|^2 + \W_2(\mu,\nu)^2\big), \ \   (x,\mu), (y, \nu)  \in \R^d\times \scr P_2.\end{split}\end{equation}  \end{enumerate}

\beg{thm}\label{T2.1} Assume {\bf (A)}.  Then  the image SDE \eqref{E1} is well-posed, and the unique solution $X_{s,t}^{x,\mu}$ is jointly continuous in $(t,x)\in [s,\infty)\times\R^d$. Moreover:
\beg{enumerate} \item[$(1)$]  For any $p\ge 1$, there  exists an increasing function $C_p: [0,\infty)\to [0,\infty)$ such that
\beq\label{EST}    \E\sup_{r\in [s,t]}\big\{ |X_{s,r}^{x,\mu}|^{2p} +\mu(|X_{s,r}^{\cdot,\mu}|^2)^{p }\big\}
 \le C_p(t)(1+|x|^{2p}+\|\mu\|_2^{2p}),\end{equation}
\beq\label{UPP}\E  \sup_{r\in [s,t]} \big\{|X_{s,r}^{x,\mu}- X_{s,r}^{y,\nu}|^{2p}+ \W_2(\LL_{s,r}^\mu, \LL_{s,r}^\nu)^{2p}\big\}
 \le C_p(t)(|x-y|^{2p}+\W_2(\mu,\nu)^{2p})\end{equation}
 hold for
all $0\le s\le t, x,y\in\R^d$ and $ \mu,\nu\in\scr P_2.$ Consequently,    $X_{s,t}^{x,\mu}$ is  jointly continuous in $(t,x)\in [s,\infty)\times \R^d$.
\item[$(2)$]   When $(b,\si)$ is deterministic,  $\{(\LL_{s,t}^\mu)_{t\ge s}: \mu\in\scr P_2\}$ is  a diffusion process on $\scr P_2$ generated  by $\BL_t$; i.e. it is a continuous strong Markov process such that for any $\mu\in\scr P_2$ and any $f\in C_b^{2}(\scr P_2)$,
$$ f(\LL_{s,t}^\mu)-f(\mu) -\int_s^t \BL_r  f(\LL_{s,r}^\mu)\d r,\ \ t\ge s$$ is a martingale.
\item[$(3)$]  When $(b,\si)$ is deterministic,  $\{(X_{s,t}^{x,\mu}, \LL_{s,t}^\mu)_{t\ge s}: \mu\in\scr P_2\}$ is  a  diffusion on $\R^d\times \scr P_2$ generated  by $\tt\BL_t$; i.e. it is a continuous strong Markov process such that for any $(x,\mu)\in\R^d\times \scr P_2$ and any $f\in C_b^{2,2}(\R^d\times \scr P_2)$,
$$ f(X_{s,t}^{x,\mu},\LL_{s,t}^\mu)-f(x,\mu) -\int_s^t \tt\BL_r f(X_{s,t}^{x,\mu},\LL_{s,r})\d r,\ \ t\ge s$$ is a martingale.
 \end{enumerate}  \end{thm}

\beg{cor} \label{C2.2} Assume {\bf (A)}.   Then for any $\F_0$-measurable random variable $X_0$ with $\mu:=\L_{X_0} \in \scr P_2$, the conditional distribution dependent SDE 
\eqref{E''} has a unique solution which is given by $X_t=X_{0,t}^{X_0,\mu}$. \end{cor}

\beg{proof} By the independence of $W$ and $\scr F_0$ and that $X_0$ is $\F_0$-measurable with distribution $\mu$,
it  is easy to show  $\L_{X_{0,t}^{X_0,\mu}|W}=\mu\circ(X_{0,t}^{\cdot,\mu})^{-1}$, which implies that   $X_t:=X_{0,t}^{X_0,\mu}$ solves \eqref{E''}.   
On the other hand, the uniqueness of \eqref{E''} can be easily proved  by using   It\^o's formula and condition \eqref{A2}.
\end{proof}

\paragraph{Remark 2.1.} Consider the conditional Mckean-Vlasov SDE (or mean field SDE with common noise)
\eqref{E'1}.  According to \cite[Theorem 3.13 and Theorem 3.35]{CD18} (see also \cite[Theorem 3.2]{CDL16} for the weak existence), if $b,\si,\si^0$ are jointly continuous on $[0,T]\times \R^d\times\scr P_2$, and there exists a constant $K>0$ such that 
$$|b(t,x,\mu)|+\|(\si,\si^0)(t,x,\mu)\|\le K(1+|x|+\|\mu\|_2),\ \ \|\nn_x b(t,x,\mu)\|+\|\nn_x(\si,\si^0)(t,x,\mu)\|\le K$$
holds for any $(t,x,\mu)\in [0,T]\times \R^d\times\scr P_2$, then for any initial distribution in $\scr P_2$, the SDE \eqref{E'1} has a weak solution up to time $T$;
if moreover $\si>0,\si^0$ are constant and $b(t,x,\mu)= b^0(t,\mu)+ cx $ for some constant $c$ and $b^0(t,\mu)$ being bounded and Lipschitz continuous in $\mu$ uniformly in $t\in [0,T]$,
then the weak solution is unique, hence the SDE is strongly well-posed. Obviously, these conditions are stronger than our assumption {\bf (A)}. 

Next, according to \cite[Theorem 2.5 and Theorem 2.7]{HSS19},  \eqref{E'1} has a weak solution provided one of the following assumptions hold: 
\beg{enumerate} \item[(i)]   $\E|X_0|^p<\infty$ for some $p>2$, $b,\si,\si^0$ are bounded and jointly continuous on $[0,T]\times \R^d\times\scr P_2$;
\item[(ii)]   $b, \si, \si^0$ are of  the integral type $$f(t,x,\mu)=\int_{\R^d} \tt f(t,x,y)\mu(\d y)$$ for bounded measurable $\tt f$, $(\tt\si,\tt\si^0) (\tt\si,\tt\si^{0})^* \ge \ll {\rm Id}$ for some constant $\ll>0$. \end{enumerate} Moreover, by \cite[Theorem 3.3]{HSS19}, if $\si$ and $\si^0$ do not depend on the distribution term such that the SDE
$$\d X_t^0= \si(t,X_t^0)\d W_t +\si^0(t,X_t^0)\d W_t^0$$ is well-posed, $\si$ is invertible such that $\{\si^{-1} b\}(t,x,\mu)$ is bounded and Lipchitz continuous in $\mu$ with respect to the total variation norm uniformly in $(t,x)\in [0,T]\times \R^d$, then the weak (hence strong) solution of \eqref{E'1} is unique.  Consequently, for the well-posedness these conditions   only apply to the non-degenerate case with  bounded  $\si^{-1}b$, but the advantage is that the drift can be   non-continuous in $(x,\mu)\in \R^d\times\scr P_2$.  

\

In the following two subsections, we prove Theorem \ref{T2.1}(1) and (2)-(3) respectively.

\subsection{Proof of Theorem \ref{T2.1}(1)}

Obviously, the uniqueness follows from \eqref{UPP}.  Below we prove  \eqref{EST},  \eqref{UPP}, joint continuity and  the existence of the solution respectively.

\paragraph{(I)  Estimate \eqref{EST}.}   Let  $(X_{s,t}^{x,\mu})_{x\in\R^d, t\ge s}$  be a solution of \eqref{E1}. We have
\beq\label{LLST} \|\LL_{s,t}^\mu\|_2^2= \|\mu\circ(X_{s,t}^{\cdot,\mu})^{-1}\|_2^2= \mu(|X_{s,t}^{\cdot,\mu}|^2),\ \ t\ge s.\end{equation} So, by \eqref{A1} and It\^o's formula, we may find out   $\kk\in L_{loc}^1([0,\infty)\to [0,\infty))$ such that
\beq\label{LLS}  \d |X_{s,t}^{x,\mu}|^2\le \kk(t) \big(1+|X_{s,t}^{x,\mu}|^2+ \mu(|X_{s,t}^{\cdot,\mu}|^2)\big)\d t + 2 \big\<X_{s,t}^{x,\mu}, \si(t,X_{s,t}^{x,\mu}, \LL_{s,t}^\mu) \d W_t\big\>,\ \ t\ge s.\end{equation}
Let $\gg_t^x= 2\si(t,X_{s,t}^{x,\mu}, \LL_{s,t}^\mu)^*X_{s,t}^{x,\mu}$.
Since $(\LL_{s,t}^\mu)_{t\ge s}$ is an adapted continuous process on $\scr P_2$ and due to \eqref{A1}, $\si(t,x,\mu)$ has linear growth in $x$,
there exists an increasing function $c: [0,\infty)\to [0,\infty)$ such that
$$\mu(|\gg_t^\cdot|)\le c(t)\big\{1+\mu(|X_{s,t}^{\cdot,\mu}|^2)\big\}= c(t)\big\{1+\|\LL_{s,t}^{\mu}\|_2^2\big\}<\infty.$$ So,    integrating  \eqref{LLS} with respect to $\mu(\d x)$ leads to
\beq\label{LLS2} \d \mu(|X_{s,t}^{\cdot,\mu}|^2)\le    \kk(t) \big(1+  2\mu(|X_{s,t}^{\cdot,\mu}|^2)\big)\d t +    \<\mu(\gg_t^\cdot),\d W_t\>, \ \ t\ge s.\end{equation}
 Let $h_{s,t} :=   \e^{2\int_s^t \kk(r)\d r}$   and
$$\tau_n= \inf\big\{t\ge s: \mu(|X_{s,t}^{\cdot,\mu}|^2) + |X_{s,t}^{x,\mu}|^2\ge n\big\},\ \ n\ge 1.$$ Then \eqref{LLS2} implies
\beq\label{OPP} \mu(|X_{s,t\land \tau_n}^{\cdot,\mu}|^2)\le h_{s,t}\|\mu\|^2_2+   \int_s^t h_{r,t}\kk(r)\d r +  \int_s^{t\land \tau_n}  h_{r,t}    \<\mu(\gg_r^\cdot),\d W_r\>, \ \ t\ge s,\end{equation} so that by \eqref{LLS},
\beq\label{PL} \beg{split} &|X_{s,t\land\tau_n}^{x,\mu}|^2\le |x|^2 + \int_s^{t\land\tau_n} \<\mu(\gg_r^\cdot), \d W_r\>\\
&+\int_s^{t\land\tau_n} \kk(r) \bigg\{1+|X_{s,r}^{x,\mu}|^2 +h_{s,r} \|\mu\|_2^2+ h_{s,r}\int_s^r \kk(\theta)\d\theta
  +\int_s^{r} h_{\theta,r} \<\gg_\theta^x,\d W_\theta\>\bigg\}\d r \end{split}\end{equation}  holds for $t\ge s$. Moreover,  \eqref{A1}  implies
\beq\label{XP} |\gg_t^x|^2  =|2\si(t,X_{s,t}^{x,\mu}, \LL_{s,t}^\mu)^*X_{s,t}^{x,\mu}|^2 \le  4 K(t) |X_{s,t}^{x,\mu}|^2\big(1+|X_{s,t}^{x,\mu}|^2+ \mu(|X_{s,t}^{\cdot,\mu}|^2)\big).\end{equation} This together with the Schwarz inequality gives
\beq\label{XP'} |\mu(\gg_t^\cdot)|^2   \le  4 K(t)\mu(|X_{s,t}^{\cdot,\mu}|^2)\big(1+2 \mu(|X_{s,t}^{\cdot,\mu}|^2)\big).\end{equation}
Then for any $p\ge 1$ and $\vv>0$, there exists a constant  $c=c(p,\vv)>0$ such that
$$\bigg(\int_s^{t\land\tau_n} | \mu(\gg_r^\cdot)|^2\d r\bigg)^{\ff p 2} \le \vv \sup_{r\in [s,t\land\tau_n]}  \big\{ \mu(|X_{s,r}^{\cdot,\mu}|^2)\big\}^p
+ c\int_s^{t\land\tau_n} K(r)\big(1+ \big\{\mu(|X_{s,r}^{\cdot,\mu}|^2)\big\}^p\big)\d r.$$
Combining this with \eqref{OPP} and using the BDG inequality, we may find  an increasing function $C_0: [0,\infty)\to\ [0,\infty)$ such that
\beg{align*}& \E\Big[\sup_{r\in [s,t\land\tau_n]}  \big\{ \mu(|X_{s,r}^{\cdot,\mu}|^2)\big\}^p\Big]\\
&\le \ff 1 2 \E\Big[\sup_{r\in [s,t\land\tau_n]}  \big\{ \mu(|X_{s,r}^{\cdot,\mu}|^2)\big\}^p\Big]+ \ff {C_0(t)} 2 \bigg(1+\|\mu\|_2^{2p}+ \E \int_s^t \big\{ \mu(|X_{s,r\land\tau_n}^{\cdot,\mu}|^2)\big\}^p\d r\bigg).\end{align*}
By Gronwall's inequality, this implies
\beq\label{OPP'} \E\Big[\sup_{r\in [s,t\land\tau_n]}  \big\{ \mu(|X_{s,r}^{\cdot,\mu}|^2)\big\}^p\Big]\le C_0(t) \e^{\int_s^t C_0(r)\d r} (1+\|\mu\|_2^{2p}).\end{equation}
Similarly, by   \eqref{PL}-\eqref{OPP'}  and  the BDG inequality,  we conclude that  for any $p\ge 1$ there exist  increasing functions $C_1,C_2: [0,\infty)\to [0,\infty)$ such that
\beg{align*} & \E \Big[\sup_{r\in [s,t\land\tau_n]} |X_{s,r}^{x,\mu}|^{2p}\Big]  \le   C_1 (t) \big(1+ |x|^{2p} +\|\mu\|_2^{2p}\big)
+ C_1(t)\E \bigg(\int_s^{t\land \tau_n} \kk(r)  |X_{s,t}^{x,\mu}|^2 \d r \bigg)^p \\
&\qquad\qquad + C_1(t)\E \bigg(\int_s^{t \land\tau_n}\big\{|\mu(\gg_r^\cdot)|^2+ |\gg_r^x|^2\big\}\d r\bigg)^{\ff p 2}\\
 &\le \ff 1 2 \E \Big[\sup_{r\in [s,t\land\tau_n]} |X_{s,r}^{x,\mu}|^{2p}\Big] + C_2(t) \big(1+ |x|^{2p} +\|\mu\|_2^{2p} \big)+ C_2(t) \E \int_s^{t }
 \kk(r)      |X_{s,r}^{x,\mu}|^{2p}  \d r ,\ \ t\ge s.\end{align*}
By Grownwall's lemma,  there exists an increasing function $Q: [0,\infty)\to (0,\infty)$ such that
$$\E \Big[\sup_{r\in [s,t\land\tau_n]} |X_{s,r}^{x,\mu}|^{2p}\Big] \le  Q(t)   (1+ |x|^{2p}+\|\mu\|_2^{2p}),\ \ t\ge s.$$
By letting $n\to \infty$ in this inequality and \eqref{OPP'},   we prove \eqref{EST} for some increasing function $C_p: [0,\infty)\to [0,\infty).$

\paragraph{(II) Estimate \eqref{UPP}.} Let $\pi\in \scr C(\mu,\nu)$ such that
\beq\label{PP1} \W_2(\mu,\nu)^2= \int_{\R^d\times\R^d} |x-y|^2\pi(\d x,\d y).\end{equation}
 Then $\pi_{s,t}:= \pi\circ (X_{s,t}^{\cdot,\mu}, X_{s,t}^{\cdot,\nu})^{-1} \in \C(\LL_{s,t}^\mu, \LL_{s,t}^\nu)$, so that
\beq\label{PP2}\beg{split} & \W_2(\LL_{s,t}^\mu,\LL_{s,t}^\nu)^2 \le \int_{\R^d\times\R^d} |x-y|^2 \pi_{s,t}(\d x,\d y)\\
&= \int_{\R^d\times\R^d} |X_{s,t}^{x,\mu}- X_{s,t}^{y,\nu}|^2 \pi(\d x,\d y)=:\ell_{s,t},\ \ t\ge s.\end{split}\end{equation}
Thus, by \eqref{A2} and It\^o's formula, we obtain
\beq\label{UO}\beg{split}  \d   |X_{s,t}^{x,\mu}- & X_{s,t}^{y,\nu}|^2   \le       K(t) \big\{|X_{s,t}^{x,\mu}-  X_{s,t}^{y,\nu}|^2+ \ell_{s,t} \big\}\d t \\
 &  + 2 \big<X_{s,t}^{x,\mu}-  X_{s,t}^{y,\nu}, \{\si(r, X_{s,t}^{x,\mu}, \LL_{s,t}^{\mu})- \si(r, X_{s,t}^{y,\nu}, \LL_{s,t}^{\nu})\}\d W_t\big\> ,\ \ t\ge  s.\end{split}\end{equation} Integrating both sides with respect to $\pi_{s,t}(\d x,\d y)$, and letting
 \beg{align*} \eta_t =  2\int_{\R^d\times\R^d}    \{\si(t, X_{s,t}^{x,\mu}, \LL_{s,t}^{\mu})- \si(t, X_{s,t}^{y,\nu}, \LL_{s,t}^{\nu})\}^*(X_{s,t}^{x,\mu}-  X_{s,t}^{y,\nu}) \,\pi(\d x,\d y),\end{align*}
 we arrive at
$$\d \ell_{s,t} \le    2K(t) \ell_{s,t}\d t + \<\eta_t,  \d W_t\>,\ \ t\ge   s.$$
This together with $\ell_{s,s}= \W_2(\mu,\nu)^2$ implies
\beq\label{IP} \ell_{s,t}  \le \W_2(\mu,\nu)^2\e^{2\int_s^tK(r)\d r} +\int_s^t \e^{2\int_r^tK(\theta)\d \theta} \<\eta_r,\d W_r\>,\ \ t\ge s.\end{equation} Moreover,     {\bf (A)} and  the Schwarz inequality yield
\beq\beg{split} \label{EAT} |\eta_r|^2&\le 4K(r)\ell_{s,r}\int_{\R^d\times \R^d} \big\{|X_{s,r}^{x,\mu}- X_{s,r}^{y,\nu}|^2+ \W_2(\LL_{s,r}^\mu, \LL_{s,r}^\nu)^2\big\}\pi(\d x,\d y)\\
&\le 8K(r) \ell_{s,r}^2,\ \ r\ge s.\end{split}\end{equation}
For given $x,y\in\R^d$ and $\mu,\nu\in \scr P_2$, let
$$\tau_n=\inf\big\{t\ge s: \|\LL_{s,t}^{\mu}\|_2+ \|\LL_{s,t}^\nu\|_2+|X_{s,t}^{x,\mu}|+ |X_{s,t}^{y,\nu}|\ge n\big\}.$$
  By  \eqref{IP}, \eqref{EAT}  and using the H\"older and BDG inequalities, we may find out     increasing functions $c_1,c_2: [0,\infty) \to [0,\infty)$ such that
 \beg{align*} &\E \Big[\sup_{r\in [s,t]} \ell_{s,r\land\tau_n}^p\Big] \le c_1(t) \W_2(\mu,\nu)^{2p} +    c_1(t) \E\bigg(\int_s^{t\land\tau_n}  |\eta_r|^2\d r\bigg)^{\ff p 2 }\\
 &\le c_1(t) \W_2(\mu,\nu)^{2p} + c_2(t) \int_s^t \E \ell_{s,r\land\tau_n}^p\d r +\ff 1 2 \E \Big[\sup_{r\in [s,t]} \ell_{s,r\land\tau_n}^p\Big],\ \ \ t\ge s.\end{align*}
Then it follows from Gronwall's lemma that
$$ \E \Big[\sup_{r\in [s,t]} \ell_{s,r\land\tau_n}^p\Big] \le 2 c_1(t) \e^{2 t c_2(t)} \W_2(\mu,\nu)^{2p},\ \ t\ge s.$$ By letting $n\to\infty$ and using Fatou's lemma, we obtain
\beq\label{ETB} \E \Big[\sup_{r\in [s,t]} \ell_{s,r}^p\Big] \le 2 c(t) \e^{2 t c_p(t)} \W_2(\mu,\nu)^{2p},\ \ t\ge s.\end{equation}
Similarly,  by  \eqref{UO},   \eqref{ETB}, assumption {\bf (A)} and using the H\"older and BDG inequality,    for any $p\ge 1$ we find out      increasing functions $K_1,K_2: [0,\infty)\to [0,\infty)$   such that
 \beg{align*}   &\E \Big[\sup_{r\in [s,t]} |X_{s,r\land\tau_n}^{x,\mu}-  X_{s,r\land\tau_n}^{y,\nu}|^{2p } \Big] \le |x-y|^{2p} +  K_1(t) \E \int_s^{t\land\tau_n} K(r) \big\{|X_{s,r}^{x,\mu}-  X_{s,r}^{y,\nu}|^{2p}+ \ell_{s,r}^{p } \big\}\d r \\
&\le |x-y|^{2p}+K_2(t) \E \int_s^{t}   K(r) |X_{s,r\land\tau_n}^{x,\mu}-  X_{s,r\land\tau_n}^{y,\nu}|^{2p}\,\d r + K_2(t)  \W_2(\mu,\nu)^{2p},
  \ \ t\ge s. \end{align*} Therefore, by Grownwall's lemma, there exists an increasing function $C: [0,\infty)\to (0,\infty)$ such that
 $$ \E\Big[\sup_{r\in [s,t]} |X_{s,r\land\tau_n}^{x,\mu}-  X_{s,r\land\tau_n}^{y,\nu}|^{2p}\Big]\le C(t) \big(|x-y|^{2p}+\W_2(\mu,\nu)^{2p}\big),\ \ t\ge s.$$
 Letting $n\to\infty$ and using Fatou's lemma, we arrive at
 $$\E\Big[\sup_{r\in [s,t]} |X_{s,r}^{x,\mu}-  X_{s,r}^{y,\nu}|^{2p}\Big]\le C(t) \big(|x-y|^{2p}+\W_2(\mu,\nu)^{2p}\big),\ \ t\ge s.$$
 Combining this  with \eqref{PP2} and \eqref{ETB}, we prove    \eqref{UPP} for some increasing function $C_p: [0,\infty)\to [0,\infty).$

\paragraph{(III) Joint continuity of $X_{s,t}^{x,\mu}$ in $(t,x).$}
 Let $K\in L_{loc}^q([0,\infty)\to [0,\infty))$ for some $q>1$. By \eqref{A1}, \eqref{EST} and \eqref{UPP}, for any $n,p\ge 1$, there exist  constants $C_1,C_2>0$ such that for any $n\ge t\ge r\ge s$, and  $|x|, |y| \le n,$
\beq\label{PPO} \beg{split} &\E\big(|X_{s,t}^{x,\mu}-X_{s,r}^{y,\mu}|^{2p} )\le 2^{2p-1}\big( \E  |X_{s,t}^{x,\mu}-X_{s,t}^{y,\mu}|^{2p} + \E |X_{s,t}^{y,\mu}-X_{s,r}^{y,\mu}|^{2p}\big)\\
&\le   C_1  |x-y|^{2p}    + C_1   \E \bigg|\int_r^t  K(\theta )\ss{1+ |X_{s,\theta}^{y,\mu}|^{2} + \mu(|X_{s,\theta}^{y,\mu}|^{2})} \d\theta\bigg|^{2p} \\
&\qquad    + C_1 \E \bigg(\int_r^t K(\theta) \big\{ 1+ |X_{s,\theta}^{y,\mu}|^{2} + \mu(|X_{s,\theta}^{y,\mu}|^{2}) \big\} \d\theta \bigg)^{p}\\
&\le   C_1  |x-y|^{2p}    + C_1  \bigg(\int_r^t K(\theta)^q \d\theta\bigg)^{\ff{2p} q} \E   \bigg|\int_r^t   \big(1+ |X_{s,\theta}^{y,\mu}|^{2}
+ \mu(|X_{s,\theta}^{y,\mu}|^{2})\big)^{\ff q{2(q-1)}} \d\theta\bigg|^{\ff{2p(q-1)}q} \\
&\qquad    + C_1\bigg(\int_r^t K(\theta)^q\bigg)^{\ff p q}  \E \bigg(\int_r^t  \big\{ 1+ |X_{s,\theta}^{y,\mu}|^{2} + \mu(|X_{s,\theta}^{y,\mu}|^{2}) \big\}^{\ff q{q-1}} \d\theta \bigg)^{\ff{p(q-1)}q}\\
 &\le  C_2\Big\{|x-y|^{2p}+  (t-r)^{\ff {p  (q-1)}q} \Big\}.  \end{split}   \end{equation} By Kolmogorov's continuity criterion, for large enough $p>1$
 this implies that $X_{s,t}^{x,\mu}$ has a $\P$-version jointly continuous in $(t,x)\in  [s,n]\times \{x\in\R^d: |x|\le n\}.$  Since $n\ge 1$ is arbitrary, $X_{s,t}^{x,\mu}$ has a version jointly continuous in $(t,x)\in  [s,\infty)\times \R^d$.

\paragraph{(IV) Existence of solution.}  It suffices to construct a solution up to an arbitrarily fixed time $T>0$.
To this end, we  adopt an iteration argument as in \cite{W18}.
\beg{enumerate} \item[(1)] For fixed $(s,\mu)\in [0,T]\times\scr P_2$, let  $\LL_{s,t}^{0,\mu}= \mu$ and $X_{s,t}^{0,x,\mu}  = x$ for all $x\in \R^d$ and $t\ge s.$
\item[(2)] Assume that for some   $n\in \mathbb Z_+$ we have constructed adapted    $(X_{s,t}^{n,x,\mu})_{t\ge s, x\in \R^d}$  which is jointly continuous in $(t,x)\in [s,\infty)\times\scr P_2$, and satisfies
\beq\label{ESTN} \E\bigg[\sup_{r\in [s,t]} |X_{s,r}^{n,x,\mu}|^2 \bigg]\le c(t)(1+|x|^2+\|\mu\|_2^2 ),\ \   t\ge s, x\in\R^d\end{equation}   for some increasing
$c: [0,\infty)\to [0,\infty)$. Consequently, $\LL_{s,t}^{n,\mu}:=
\mu\circ(X_{s,t}^{n,\cdot,\mu})^{-1}\in \scr P_2$ is   continuous in $t\ge s$. Indeed, by the Fubini theorem, \eqref{ESTN} implies
$$\E\bigg[\mu\Big(\sup_{r\in [s,t]} |X_{s,r}^{n,\cdot,\mu}|^2 \Big)\bigg]\le c(t)(1+ 2\|\mu\|_2^2)<\infty,\ \ t\ge s,$$ so that $\P$-a.s
$$\mu\Big(\sup_{r\in [s,t]} |X_{s,r}^{n,\cdot,\mu}|^2 \Big)<\infty,\ \ t\ge s.$$
Then by the dominated convergence theorem and the continuity of $X_{s,t}^{n,x,\mu}$ in $t\ge s$, we obtain $\P$-a.s.
$$\lim_{r\to t} \W_2(\LL_{s,r\lor s}^{n,\mu}, \LL_{s,t}^{n,\mu})^2\le \lim_{r\to t} \mu\big(|X_{s,r\lor s}^{n,\cdot,\mu}- X_{s,t}^{n,\cdot,\mu}|^2\big)=0,\ \ t\ge s.$$
\item[(3)] Let
$(X_{s,t}^{n+1, x,\mu})_{t\ge s} $ solve the SDE
$$\d X_{s,t}^{n+1,x,\mu}= b(t,X_{s,t}^{n+1,x,\mu},  \LL_{s,t}^{n,\mu})\d t +\si(t,X_{s,t}^{n+1,x,\mu},\LL_{s,t}^{n,\mu})\d W_t,\ \ t\ge s, X_{s,s}^{n+1,x,\mu}=x.$$
By {\bf (A)} and \eqref{ESTN}, it is easy to see that   this SDE  is well-posed, and when $x$ varies the inequality \eqref{ESTN}  holds for  $X_{s,t}^{n+1,x,\mu}$  replacing $X_{s,t}^{n,x,\mu}$ with possibly a different function $c: [0,\infty)\to [0,\infty)$. Moreover, as  in {\bf (III)}, {\bf (A)} and \eqref{ESTN} also imply the joint continuity of $X_{s,t}^{n+1,x,\mu}$ in $(t,x)\in [s,\infty)\times\R^d$. Consequently, as shown in step (2) that
$\LL_{s,t}^{n+1,\mu}:=\mu\circ(X_{s,t}^{n+1,\cdot,\mu})^{-1}\in\scr P_2$ is continuous in $t\ge s$.   \end{enumerate}

Therefore, we have constructed a sequence $\{(X_{s,t}^{n,x,\mu}, \LL_{s,t}^{n,\mu})_{t\ge s, x\in\R^d}\}_{n\ge 0}$, which satisfies \eqref{ESTN},
$X_{s,t}^{n,x,\mu}$ is jointly continuous in $(t,x)\in [s,\infty)\times \R^d$, and $\P$-a.s.
\beq\label{AL} X_{s,t}^{n+1,x,\mu}= x+\int_s^t b(r,X_{s,r}^{n+1,x,\mu},  \LL_{s,r}^{n,\mu})\d r +\int_s^t\si(r,X_{s,r}^{n+1,x,\mu},\LL_{s,r}^{n,\mu})\d W_r,\ \ t\ge s, x\in\R^d.\end{equation}
The following lemma gives a constant $t_0>0$ independent of $(s,x,\mu)\in [0,T]\times\R^d\times\scr P_2$,  such that   $\{X_{s,\cdot}^{n,x,\mu}\}_{n\ge 1}$ is a Cauchy sequence in $L^2(\OO\to C([s,s+t_0]\to\R^d); \P)$.

\beg{lem} \label{LL1} Assume {\bf (A)}.  For fixed $T>0,$ there exists a  constant  $t_0>0$ such that
 $$\lim_{n,m\to\infty} \sup_{(s,x,\mu)\in [0,T]\times\R^d\times\scr P_2} \ff{\E \sup_{t\in [s,s+t_0]} |X_{s,t}^{m,x,\mu}-X_{s,t}^{n,x,\mu}|^2}{1+|x|^2+\|\mu\|_2^2}   =0.$$
\end{lem}

\beg{proof} As in \eqref{PP2}, we have $\W_2(\LL_{s,t}^{n,\mu}, \LL_{s,t}^{n-1,\mu})^2 \le \mu(|X_{s,t}^{n,\cdot,\mu} - X_{s,t}^{n-1,\cdot,\mu}|^2)$ for $n\ge 1$. Combining this with  \eqref{A2}  and It\^o's formula, we obtain \beg{align*}& \d |X_{s,t}^{n+1,x,\mu}-  X_{s,t}^{n,x,\mu}|^2 \le   K(t) \Big\{|X_{s,t}^{n+1,x,\mu}-  X_{s,t}^{n, x,\mu}|^2+ \mu(|X_{s,t}^{n,\cdot,\mu}-   X_{s,t}^{n-1,\cdot,\mu}|^2)\Big\}\d t\\
&\qquad  + 2\big<X_{s,t}^{n+1,x,\mu}-  X_{s,t}^{n, x,\mu}, \{\si(t, X_{s,t}^{n+1,x,\mu}, \LL_{s,t}^{n,\mu})- \si(t, X_{s,t}^{n,x,\mu}, \LL_{s,t}^{n-1,\mu})\}\d W_t\big\> ,\ \ t\ge s.\end{align*}  So, by \eqref{A2} and using the BDG inequality,  we may find out   constants $c_1,c_2>0$   such that
\beg{align*}  &\E \Big[\sup_{t\in [s, s+t_0]} |X_{s,t}^{n+1,x,\mu}-  X_{s,t}^{n,x,\mu}|^2 \Big]\\
&\le   \int_s^t K(r) \E \big[|X_{s,r}^{n+1,x,\mu}-  X_{s,r}^{n,x,\mu}|^2+ \mu(|X_{s,r}^{n,\cdot,\mu}-  X_{s,t}^{n-1,\cdot,\mu}|^2)\big] \d r\\
& \quad + c_1\E\bigg(\int_s^t K(r) |X_{s,r}^{n+1,x,\mu}-  X_{s,r}^{n,x,\mu}|^2\big\{|X_{s,r}^{n+1,x,\mu}-  X_{s,t}^{n,x,\mu}|^2+\mu(|X_{s,r}^{n,\cdot,\mu}-  X_{s,t}^{n-1,\cdot,\mu}|^2)\big\} \d r\bigg)^{\ff 1 2}\\
&\le \ff {c_2}  2    \int_s^t K(r) \E \big[|X_{s,r}^{n+1,x,\mu}-  X_{s,r}^{n,x,\mu}|^2+ \mu(|X_{s,r}^{n,\cdot,\mu}-  X_{s,t}^{n-1,\cdot,\mu}|^2)\big] \d r\\
&\quad  +\ff 1 2 \E \Big[\sup_{t\in [s, s+t_0]} |X_{s,t}^{n+1,x,\mu}-  X_{s,t}^{n,x,\mu}|^2 \Big],\ \ t\ge s.\end{align*}
Since \eqref{ESTN}  holds for all $n$, this and Grownwall's inequality imply
\beq\label{TTT} \E \sup_{r\in [s, t]} |X_{s,r}^{n+1,x,\mu}-  X_{s,r}^{n,x,\mu}|^2\le c_2\int_s^{t}\e^{c_2\int_r^{t} K(\theta)\d\theta} \E\mu(|X_{s,r}^{n,\cdot,\mu}-  X_{s,r}^{n-1,\cdot,\mu}|^2)\d r,\ \ t\ge s\end{equation}    for all $(s,x)\in [0,T]\times \R^d$. Taking integral with respect to $\mu(\d x)$ leads to
$$\sup_{r\in [s,t]} \E\mu( |X_{s,r}^{n+1,\cdot,\mu}-  X_{s,r}^{n,\cdot,\mu}|^2)\le c_2 (t-s)\e^{c_2\int_s^{t}K(r)\d r} \sup_{r\in [s,t]} \E\mu( |X_{s,r}^{n,\cdot,\mu}-  X_{s,r}^{n-1,\cdot,\mu}|^2),\ \ t\ge s.$$
Now, taking $t_0>0$ such that
\beq\label{TVV} \vv:= c_2 t_0 \e^{c_2\int_0^{T+t_0} K(r)\d r}<1,\end{equation}
by iterating in $n$ we arrive at
\beg{align*} & \sup_{s\in [0,T], t\in [s,s+t_0]} \E\mu( |X_{s,t}^{n+1,\cdot,\mu}-  X_{s,t}^{n,\cdot,\mu}|^2)\le \vv \sup_{s\in [0,T], t\in [s,s+t_0]} \E\mu( |X_{s,t}^{n,\cdot,\mu}-  X_{s,t}^{n-1,\cdot,\mu}|^2) \\
&\le\cdots\le  \vv^n \sup_{s\in [0,T], t\in [s,s+t_0]}
\E\mu( |X_{s,t}^{1,\cdot,\mu}-  X_{s,t}^{0,\cdot,\mu}|^2)= c(x,\mu)\vv^n<\infty,\end{align*}
where due to \eqref{ESTN},
$$c(x,\mu):= \sup_{s\in [0,T]}\sup_{t\in [s,s+t_0]}
\E\mu( |X_{s,t}^{1,\cdot,\mu}-  x|^2)\le c(1+|x|^2+\|\mu\|_2^2)$$ for some  constant $c>0$.
Substituting this into \eqref{TTT} and using \eqref{TVV}, we get
$$ \sup_{s\in [0,T]} \E \sup_{t\in [s, s+t_0]} |X_{s,t}^{n+1,x,\mu}-  X_{s,t}^{n,x,\mu}|^2\le c(1+|x|^2+\|\mu\|_2^2) \vv^n,\ \ n\ge 1.$$
This finishes the proof.
\end{proof}

By Lemma \ref{LL1}, there exist  a constant $t_0>0$ depending on $T>0$, such that for any $s\in [0,T)$ we have a family of continuous processes
$$\{(X_{s,t}^{x,\mu})_{t\in [s,s+t_0]}: x\in\R^d, \mu\in\scr P_2\}$$ which are measurable in $x$  and
 $$\lim_{n\to\infty} \E\Big[\sup_{r\in [s,s+t_0]} \big(|X_{s,r}^{n,x,\mu} - X_{s,t}^{x,\mu}|^2 + \mu(|X_{s,r}^{n,\cdot,\mu} - X_{s,t}^{\cdot,\mu}|^2)\big)\Big]=0.$$
Letting $\LL_{s,t}^\mu= \mu\circ(X_{s,t}^{\cdot,\mu})^{-1}$, by this and    \eqref{PP2} we obtain
   $$\lim_{n\to\infty} \E\Big[\sup_{r\in [s,s+t_0]} \W_2(\LL_{s,t}^{n,\mu}, \LL_{s,t}^\mu)^2\Big]\le \E \Big[\sup_{r\in [s,s+t_0]}  \mu(|X_{s,r}^{n,\cdot,\mu} - X_{s,t}^{\cdot,\mu}|^2) \Big]=0.$$  Thus, the continuity of $\LL_{s,t}^{n,\mu}$ in $t\in [s, s+t_0]$ implies that of $\LL_{s,t}^\mu$; due to   \eqref{ESTN}  we may find out a   constant $c_1>0$ such that
\beq\label{AL2} \E \Big[\sup_{t\in [s,s+t_0]} \big\{\mu(|X_{s,t}^{\cdot,\mu}|^2)+ |X_{s,t}^{x,\mu}|^2\big\}\Big] \le c_1 \big(1+|x|^2+\|\mu\|_2^2\big),\ \ (s,x,\mu)\in [0,T]\times  \R^d\times\scr P_2;\end{equation}  and finally,  by  assumption {\bf (A)} we may let $n\to\infty$ in \eqref{AL} to derive
$$X_{s,t}^{x,\mu} = x + \int_s^t b(r,X_{s,r}^{x,\mu},\LL_{s,r}^\mu)\d r+ \int_s^t \si(r, X_{s,r}^{x,\mu},\LL_{s,r}^\mu)\d W_r,\ \ t\in [s, s+t_0], x\in \R^d.$$
 So, when $T\le s+t_0$ we have   solved the SDE up to time $T.$

 In the case that $T>s+t_0$,   let $\bar s= s+t_0, \bar x= X_{s,s+t_0}^{x,\mu}$ and $ \bar\mu= \LL_{s,s+t_0}^\mu.$ Since
given $\F_{s+t_0}$ the process $(W_t-W_{\bar s})_{t\ge \bar s}$ is an $m$-dimensional Brownian motion, and $(\bar x,\bar\mu)$ is given as well, as in above  we may construct a solution $(X_{\bar s, t}^{\bar x,\bar \mu}, \LL_{\bar s, t}^{\bar\mu})_{t\in [\bar s, \bar s+t_0]} $ for \eqref{E1} with $\bar s$ replacing $s$. Then extending $(X_{s,t}^{x,\mu}, \LL_{s,t}^\mu)$ to $t\in [\bar s,\bar s+t_0]$ by letting
$$X_{s,t}^{x,\mu}= X_{\bar s, t}^{\bar x, \bar\mu},\ \ \LL_{s,t}^\mu= \LL_{\bar s, t}^{\bar \mu},\ \ t\in [\bar s, \bar s+t_0],$$
we see that $(X_{s,t}^{x,\mu}, \LL_{s,t}^\mu)_{t\in [s, s+2t_0]}$ solves \eqref{E1} up to time $\bar s+t_0=s+2t_0$. Runing this procedure for $k$ times until $s+kt_0\ge T$, we construct a solution to \eqref{E1} up to time $T$.

\subsection{Proof of Theorem \ref{T2.1}(2)-(3)}

We first establish   It\^o's formula for the diffusion process $(\LL_{s,t}^\mu)_{t\ge s}$. To this end, we need the following chain rule for the $L$-derivative, which is essentially due to \cite[Theorem 6.5]{Card} where the reference probability space is Polish,   see also \cite[Proposition A.2]{HSS} for general probability space but bounded random variables $\{\xi_s\}_{s\in [0,\vv]}$ (note that $D_k$ therein is compact).

\beg{lem}\label{L2.0} Let $\{\xi_s\}_{s\in [0,\vv]}$ for some $\vv>0$ be a family of square integrable random variables on $\R^d$ with respect to a probability space $(\OO,\scr F,\P)$, and let $\L_{\xi_s}$ denote the law of $\xi_s$. If
$$ \xi_0':=\lim_{s\downarrow 0} \ff{\xi_s-\xi_0}s$$ exists in $L^2(\OO\to\R^d;\P)$, then for any   $f\in C^{1}(\scr P_2)$,
$$\lim_{s\downarrow 0} \ff{f(\L_{\xi_s})- f(\L_{\xi_0})}s = \E\<Df(\L_{\xi_0})(\xi_0), \xi_0'\>.$$\end{lem}

\beg{proof} By a standard extension argument, we may and do assume that $(\OO,\scr F,\P)$ is atomless.  For instance, we enlarge   $(\OO,\scr F,\P)$ by $(\OO\times [0,1],\scr F\times\scr B([0,1]),\P\times\d r)$ and use $\tt \xi_s$ to replace $\xi_s$, where $\tt\xi_s(\oo,r):=\xi_s(\oo)$ for $(\oo,r)\in \OO\times [0,1],$ so that $\L_{\tt\xi_s}$ under $\P\times\d r$ coincides with $\L_{\xi_s}$ under $\P$.  Then the  proof is completely similar to that of \cite[Proposition 3.1]{RW18b} for $\xi_s$ replacing $X+sY$.
\end{proof}

\beg{lem}[It\^o's formula] \label{L2.1} Assume {\bf (A)} and let $\{\LL_{s,t}^\mu= \mu\circ(X_{s,t}^{\cdot,\mu})^{-1}\}_{t\ge s}$ for  the  solution to $\eqref{E1}$.       Then   for any $f\in C_b^{2}(\scr P_2)$,
$$\d f(\LL_{s,t}^\mu)= (\BL_t  f)(\LL_{s,t}^\mu) \d t + \bigg\<\int_{\R^d}\big\{\si(t,x,\LL_{s,t}^\mu)^* (Df)(\LL_{s,t}^\mu)(x)\big\} \mu(\d x), \d W_t\bigg\>,\ \ t\ge s.$$
\end{lem}

\beg{proof} For any $t\ge s$ and small $\vv>0$, let
$$\xi_r=  (1-r)X_{s,t}^{\cdot,\mu} + r X_{s,t+\vv}^{\cdot,\mu}: \R^d\to\R^d,\ \ r\in [0,1].$$
Then $\mu\circ\xi_r^{-1}$ is the law of $\xi_r$ on the probability space $(\R^d,\B(\R^d),\mu).$ By \eqref{EST},
$$\sup_{r\in [0,1]} \E \|\mu\circ\xi_r^{-1}\|_2^2 \le     \E\Big[\sup_{r\in [0,1]}  \mu(|\xi_r|^2) \Big] <\infty,\ \ t\ge s.$$
 Moreover,
$\xi_r':=\ff{\d}{\d r}\xi_r= X_{s,t+\vv}^{\cdot,\mu}- X_{s,t}^{\cdot,\mu}$ exists in $L^2(\R^d\to\R^d;\mu)$. So,  Lemma \ref{L2.0}  implies
\beq\label{BG1} \beg{split} &f(\LL_{s,t+\vv}^\mu)- f(\LL_{s,t}^\mu)  =  f(\mu\circ\xi_1^{-1} )- f(\mu\circ \xi_0^{-1})
= \int_0^1 \Big(\ff{\d}{\d r}   f(\mu\circ\xi_r^{-1} ) \Big) \d r\\
&= \int_{\R^d\times [0,1]} \big\<Df(\mu\circ\xi_r^{-1})(\xi_r^x), X_{s,t+\vv}^{x,\mu}-X_{s,t}^{x,\mu}\big\>\mu(\d x)\d r\\
&= \int_{\R^d} I_1(x)\mu(\d x)+\int_{\R^d\times[0,1]} I_2(x,r)\mu(\d x)\d r +\int_{\R^d\times [0,1]} I_3(x,r)\mu(\d x)\d r,
\end{split} \end{equation} where, since $\mu\circ\xi_0^{-1}=\LL_{s,t}^\mu$,
\beg{align*} &I_1(x):=   \<D f(\LL_{s,t}^\mu)(X_{s,t}^{x,\mu}), X_{s,t+\vv}^{x,\mu}- X_{s,t}^{x,\mu}\big\>,\\
&I_2(x,r):=   \big\<D f(\mu\circ \xi_r^{-1})(\xi_r^x)- Df(\mu\circ\xi_0^{-1}) (\xi_r^x), X_{s,t+\vv}^{x,\mu}- X_{s,t}^{x,\mu}\big\>,  \\
&I_3(x,r):=   \big\<D f(\LL_{s,t}^\mu)(\xi_r^x)- Df(\LL_{s,t}^\mu) (\xi_0^x), X_{s,t+\vv}^{x,\mu}- X_{s,t}^{x,\mu}\big\>.
\end{align*}
Below, we calculate $I_1(x), I_2(x)$ and $I_3(x)$ respectively.

Firstly,
by \eqref{E1} and $f\in C_b^2(\scr P_2)$, we have
\beq\label{BG2} \beg{split} I_1 (x) &= \int_t^{t+\vv} \big\<(D f)(\LL_{s,t}^\mu)(X_{s,t}^{x,\mu}), \d X_{s,u}^{x,\mu}\big\> = \int_t^{t+\vv} \big\<(D f)(\LL_{s,u}^\mu)(X_{s,u}^{x,\mu}), \d X_{s,u}^{x,\mu}\big\> + {\rm o}(\vv)\\
&= \int_t^{t+\vv} \big\<(D f)(\LL_{s,u}^\mu)(X_{s,u}^{x,\mu}), b(u, X_{s,u}^{x,\mu},\LL_{s,u}^\mu)\big\>\d u \\
&\qquad +  \int_t^{t+\vv} \big\<(D f)(\LL_{s,u}^\mu)(X_{s,u}^{x,\mu}), \si (u,\d X_{s,u}^{x,\mu},\LL_{s,u}^\mu)\d W_u\big\>+
   {\rm o}(\vv),\end{split} \end{equation}where and in the following, ${\rm o}(\vv)$ means $\vv$-dependent (real, vector or matrix valued) random variables satisfying
$\lim_{\vv\to 0} \vv^{-1}   |{\rm o}(\vv)|=0.$

Next, \eqref{E1} implies
\beq\label{**0} (X_{s,t+\vv}^{x,\mu}- X_{s,t}^{x,\mu})\otimes (X_{s,t+\vv}^{y,\mu}- X_{s,t}^{y,\mu})= \int_t^{t+\vv}   \si(u, X_{s,u}^{x,\mu}, \LL_{s,u}^\mu)
\si(u,X_{s,u}^{y,\mu},\LL_{s,u}^\mu)^*\d u +  {\rm o}(\vv).\end{equation}
Combining this with $f\in C_b^2(\scr P_2)$, we deduce from   Lemma \ref{L2.0}     and $\xi_\theta'=X_{s,t+\vv}^{\cdot,\mu}- X_{s,t}^{\cdot,\mu}$  that up to an error term ${\rm o}(\vv)$,
\beq\label{BG3} \beg{split} &I_2 (x,r)= \int_0^r \d\theta \int_{\R^d} \big\<(D^2 f)(\mu\circ \xi_\theta^{-1})(\xi_r^x, \xi_\theta^y), (X_{s,t+\vv}^{x,\mu}- X_{s,t}^{x,\mu})\otimes (\xi_\theta^y)'\big\>\mu(\d y) \\
&= \int_0^r \d\theta\int_t^{t+\vv} \d u \int_{\R^d} \big\<(D^2) f(\mu\circ \xi_\theta^{-1})(\xi_r^x, \xi_\theta^y), \si(u, X_{s,u}^{x,\mu}, \LL_{s,u}^\mu)\si(u,X_{s,u}^{y,\mu},\LL_{s,u}^\mu)^*\big\>\mu(\d y) \\
&=  r \int_t^{t+\vv}\d u \int_{\R^d}  \big\<(D^2 f)(\LL_{s,u}^\mu)(X_{s,u}^{x,\mu}, X_{s,u}^{y,\mu}),
\si(u, X_{s,u}^{x,\mu}, \LL_{s,u}^\mu)\si(u,X_{s,u}^{y,\mu},\LL_{s,u}^\mu)^*\big\>\mu(\d y).\end{split}  \end{equation}

Similarly,    by using \eqref{**0} with $x=y$, we obtain that  up to an error term ${\rm o}(\vv)$,
\beg{align*} I_3(x,r)&=  \big\<(D f)(\mu\circ \xi_0^{-1})(\xi_r^x)- (Df)(\mu\circ\xi_0^{-1}) (\xi_0^x), X_{s,t+\vv}^{x,\mu}- X_{s,t}^{x,\mu}\big\>\\
&=\int_0^r \big\<\nn \{(D f)(\LL_{s,t}^\mu)\}(\xi_\theta^x), (X_{s,t+\vv}^{x,\mu}- X_{s,t}^{x,\mu})\otimes(X_{s,t+\vv}^{x,\mu}- X_{s,t}^{x,\mu})\big\> \d r\\
&=  r\int_t^{t+\vv}    \big\<\nn \{(D f)(\LL_{s,u}^\mu)\}(X_{s,u}^{x,\mu}),   (\si\si^*)(t, X_{s,u}^{x,\mu}, \LL_{s,u}^\mu)\big\>\d u.\end{align*}
Combining this with \eqref{BG1}-\eqref{BG3}, we arrive at
\beg{align*} &\d f(\LL_{s,t}^\mu)- \int_{\R^d}\big\< (Df)(\LL_{s,t}^\mu)(x), \si(t,x,\LL_{s,t}^\mu)\d W_t\big\>\mu(\d x)\\
&=  \bigg(\int_{\R^d}   \big\<(D f)(\LL_{s,t}^\mu)(X_{s,t}^{x,\mu}), b(t, X_{s,t}^{x,\mu}, \LL_{s,t}^\mu)\big\>\mu(\d x) \bigg)\d t \\
&\qquad +
   \bigg(\ff 1 2\int_{\R^d} \big\<\nn \{(D f)(\LL_{s,t}^\mu)\}(X_{s,t}^{x,\mu}),   (\si\si^*)(t, X_{s,t}^{x,\mu}, \LL_{s,t}^\mu)\big\>\Big\}\mu(\d x) \bigg)\d t\\
  &\qquad + \bigg(\ff 1 2 \int_{\R^d\times\R^d}  \big\<(D^2 f)(\LL_{s,t}^\mu)(X_{s,t}^{x,\mu}, X_{s,t}^{y,\mu}),
\si(t, X_{s,t}^{x,\mu}, \LL_{s,t}^\mu)\si(t,X_{s,t}^{y,\mu},\LL_{s,t}^\mu)^*\big\>\mu(\d x)\mu(\d y)\bigg)\d t\\
&= (\BL_t  f)(\LL_{s,t}^\mu) \d t.\end{align*} Then the proof is finished.
 \end{proof}
 
 \paragraph{Remark 2.2.}  We note that under  a moment condition the the I\^o's formula for the conditional distribution has been established in \cite{CD18}. More precisely, consider the It\^o process 
 $$\d X_t= B_t\d t+ \Sigma_t\d W_t+\Sigma^0_t\d W_t^0,$$ where 
 $B_t, \Sigma_t$ and $\Sigma_t^0$ are progressively measurable such that 
 \beq\label{4M} \E\int_0^T \big\{|B_t|^2+ \|\Sigma_t\|^4+ \|\Sigma_t^0\|^4\big\}\d t<\infty,\end{equation}
 then for any $f\in C^2_b(\scr P_2)$, $\mu_t:=\L_{X_t|W}$ satisfies the It\^o formula
 \beg{align*} & f(\mu_t)-f(\mu_0) = \int_0^t \E\big[\<Df(\mu_s)(X_s),B_s\>|W\big]\d s +\int_0^t \big\<\E\big[Df(\mu_s)(X_s)\big|W\big], \Sigma_s\d W_s\big\>\\
 &+\ff 1 2 \int_0^t \E\big[{\rm tr}\big\{\nn \big(Df(\mu_s)\big)(X_s)\big(\Sigma_s\Sigma_s^* +\Sigma_s^0 (\Sigma_s^0)^*\big)\big|W\big]\d s \\
 &+\ff 1 2 \int_0^t \d s \int_{\R^d} \E\big[{\rm tr}\big\{D^2f(\mu_s)(X_s,z) \Sigma_s\Sigma_s^*\big\}\big|W\big] \mu_s(\d z),\ \ t\in [0,T].\end{align*}
So, if \eqref{4M} holds for $B_t:= b(t, X_t,\L_{X_t|W}), \Sigma_t:= \si(t,X_t, \L_{X_t|W})$ and $\Sigma_t^0=0$, then Lemma \ref{L2.1} follows. However, our condition {\bf (A)} is not enough to ensure \eqref{4M} unless the initial value satisfies $\E|X_0|^4<\infty$ (as $\si$ has linear growth) and $K\in L^2([0,T])$.

\beg{proof}[Proof of Theorem \ref{T2.1}(2)-(3)]     By the uniqueness result  in Theorem \ref{T2.1}, we have the flow property
\beq\label{MK}X_{s,t}^{x,\mu}= X_{r,t}^{X_{s,r}^{x,\mu}, \LL_{s,r}^\mu},\ \ \LL_{s,t}^\mu= \LL_{r,t}^{\LL_{s,r}^\mu},\ \ 0\le s\le r \le t,\end{equation}
which implies that both $(\LL_{s,t}^\mu)_{t\ge s} $ and $(X_{s,t}^{x,\mu}, \LL_{s,t}^\mu)_{t\ge s}$ are Markov processes.

Next, by \eqref{UPP},
these two Markov processes are Feller and hence,   strong Markovian.       Therefore, Theorem \ref{T2.1}(2) follows from Lemma \ref{L2.1}.

Finally, for any $f\in C_b^{2,2}(\R^d,\scr P_2)$,
Lemma \ref{L2.1} and the classical It\^o's formula for the semimartingale $(X_{s,t}^{x,\mu})_{t\ge s}$ imply
\beq\label{ITO}   \beg{split}    \d f(X_{s,t}^{x,\mu}, \LL_{s,t}^\mu)   = &\ (\tt\BL_t f)(X_{s,t}^{x,\mu}, \LL_{s,t}^\mu)  \d t  + \big\<\nn f(\cdot,\LL_{s,t}^{\mu})(X_{s,t}^{x,\mu}), \si(t,X_{s,t}^{x,\mu}, \LL_{s,t}^{\mu})\d W_t\big\>\\
&\  + \int_{\R^d}\big\< Df(X_{s,t}^{x,\mu},\cdot)(\LL_{s,t}^\mu)(x), \si(t,x,\LL_{s,t}^\mu)\d W_t\big\>\mu(\d x),    \ \ t\ge s.\end{split} \end{equation}
This proves Theorem \ref{T2.1}(3).\end{proof}

\section{Feynman-Kac formula for PDEs on $\R^d\times \scr P_2$}

In this section, we solve the   PDEs \eqref{PDE} and \eqref{PDE2} by using  $(X_{s,t}^{x,\mu},\LL_{s,t}^\mu)_{0\le s\le t\le T}.$
As mentioned in Abstract that when $V=0$ they  are included by the Master equations investigated in the literature of mean filed games with common noise. 

A function on $U$ on $[0,T]\times\R^d\times\scr P_2$ is called a solution to \eqref{PDE}, if $U(t,x,\mu)$ is differentiable in $t$ and $U(t,\cdot,\cdot)\in C^{2,2}(\R^d\times\scr P_2)$ such that \eqref{PDE} holds. If moreover $U(t,x,\mu)$ does not depend on $x$, it is called a solution to \eqref{PDE2}.
We first introduce the following  class  $C_b^{0,2,2}([0,T] \times\R^d\times \scr P_2)$.

\beg{defn}   Let $f$ be a real, vector or matrix valued function on
$[0,T]\times \R^k\times \scr P_2$ for some $k\ge 1$.   We write   $f\in C_b^{0,2,2}([0,T] \times\R^k\times \scr P_2)$, if $f$ is jointly continuous, $f(t,\cdot,\cdot)\in C_b^{2,2}(\R^k\times \scr P_2)$ for every $t\in [0,T]$, and all derivatives
\beg{align*}& \nn f(t,x,\mu),\ \ \nn^2 f(t,x,\mu),\ \ D f(t,x,\mu)(y),\\
& D\{\nn f(t,x,\mu)\}(y),\ \ \nn \{Df(t,x,\mu)(\cdot)\}(y),\ \ D^2f(t,x,\mu)(y,z)\end{align*}
are bounded and jointly continuous in corresponding arguments. If moreover  $f(t,x,\mu)$ does not depend on $x$, we denote $f\in C_b^{0,2}([0,T]\times \scr P_2).$
\end{defn}

\beg{thm}\label{T4.1} Assume that $b,\si\in C_b^{0,2,2}([0,T] \times\R^d\times \scr P_2)$ are deterministic.
\beg{enumerate} \item[$(1)$] For any $\Phi\in C_b^{2,2}(\R^d\times\scr P_2)$, $F\in  C_b^{0,2,2}([0,T] \times\R^d\times \scr P_2)$, and bounded
$V\in C_b^{0,2,2}([0,T] \times\R^d\times \scr P_2)$,
$$U(t,x,\mu):= \E \bigg[\Phi(X_{t,T}^{x,\mu}, \LL_{t,T}^\mu)\e^{\int_t^T V(r,X_{t,r}^{x,\mu},\LL_{t,r}^\mu)d r}+
\int_t^T F(r,X_{t,r}^{x,\mu}, \LL_{t,r}^\mu)\e^{\int_t^r V(\theta, X_{t,\theta}^{x,\mu},\LL_{t,\theta}^\mu)\d\theta}\d r \bigg]$$ is the unique solution of \eqref{PDE} in the class $C_b^{0,2,2}([0,T] \times\R^d\times \scr P_2)$ with $\pp_t U\in C([0,T]\times\R^d\times\scr P_2)$.
\item[$(2)$] For any $\Phi\in C_b^{2}(\R^d\times\scr P_2)$, $ F\in  C_b^{0,2}([0,T]  \times \scr P_2)$, and bounded
$V\in C_b^{0,2}([0,T] \times \scr P_2)$,
$$U(t,\mu):= \E \Big[\Phi(\LL_{t,T}^\mu)\e^{\int_t^T V(r,\LL_{t,r}^\mu)\d r}+
\int_t^T  F(r, \LL_{t,r}^\mu)\e^{\int_t^r V(\theta, \LL_{t,\theta}^\mu)\d\theta}\d r \Big]$$ is the unique solution of \eqref{PDE2} in the class $C_b^{0,2}([0,T] \times \scr P_2)$   with $\pp_t U \in C([0,T]\times\scr P_2)$ \end{enumerate}
\end{thm}

\paragraph{Remark 3.1.}  When $\si$ is constant, $b(t,x,\mu)$ and $\nn_x b(t,x,\mu)$ are in the class $C_b^{0,2,2}([0,T]\times\R^d\times\scr P_2)$ and $V=0$, 
\cite[Theorem 5.45]{CD18} implies that $U(t,x,\mu)$ given in Theorem \ref{T4.1}(1) solves  the Master equation \eqref{PDE} with $V=0$. These conditions are stronger than those in Theorem \ref{T4.1}.  

\beg{proof}[Proof of Theorem \ref{T4.1}] Since   $\tt \BL_t  F(x,\mu) = \BL_t  F(\mu)$ holds for $F\in C_b^{2}(\scr P_2)$,
(2) follows from (1).
So, it suffices to prove Theorem \ref{T4.1}(1).

If $U\in C_b^{0,2,2}([0,T] \times\R^d\times \scr P_2)$ is a solution of \eqref{PDE}, then  \eqref{ITO}  yields
\beg{align*} \d U(t, X_{s,t}^{x,\mu}, \LL_{s,t}^\mu)& = (\pp_t + \tt \BL_t) U(t, X_{s,t}^{x,\mu}, \LL_{s,t}^\mu) \d t + \d M_t\\
&= \d M_t-\big(VU+F\big)(t, X_{s,t}^{x,\mu}, \LL_{s,t}^\mu)\d t, \ \ t\in [s,T] \end{align*}  for some  martingale $(M_t)_{t\in [s,T]}.$  Thus,   the process $$\eta_t:=
  U(t, X_{s,t}^{x,\mu}, \LL_{s,t}^\mu) \e^{\int_s^t V(r,X_{s,r}^{x,\mu}, \LL_{s,r}^\mu)\d r} + \int_s^t F(r,X_{s,r}^{x,\mu}, \LL_{s,r}^\mu)\e^{\int_s^r V(\theta,X_{s,\theta}^{x,\mu}, \LL_{s,\theta}^\mu)\d\theta}\d r,\ \ t\in [s,T]$$ satisfies
  $$\d \eta_t= \e^{\int_s^t V(r,X_{s,r}^{x,\mu}, \LL_{s,r}^\mu)\d r}\d M_t,\ \ t\in [s,T].$$
So, \beg{align*} &U(s,x,\mu)= \E\eta_s = \E \eta_T  \\
&= \E \bigg[\Phi(X_{s,T}^{x,\mu}, \LL_{s,T}^\mu)\e^{\int_s^T V(r,X_{s,r}^{x,\mu},\LL_{s,r}^\mu)d r}+
\int_s^T F(r,X_{s,r}^{x,\mu}, \LL_{s,r}^\mu)\e^{\int_s^r V(\theta, X_{s,\theta}^{x,\mu},\LL_{s,\theta}^\mu)\d\theta}\d r \bigg]  \end{align*}
 as claimed in Theorem \ref{T4.1}(1).

On the other hand, let $U$ be given in Theorem \ref{T4.1}(1). For any  $t\in [0,T)$ and $\vv\in (0,T-t)$, by  \eqref{MK} and the formula of $U(t,x,\mu)$ in Theorem \ref{T4.1}(1),
\beg{align*} U(t,x,\mu)-\E \big[U(t+\vv, X_{t,t+\vv}^{x,\mu},\LL_{t,t+\vv}^x) \big]
 =  I_1(\vv)+I_2(\vv)+I_3(\vv) \end{align*} holds for
  \beg{align*} &I_1(\vv):=   \E \Big[\Phi(X_{t+\vv,T}^{X_{t,t+\vv}^{x,\mu}, \LL_{t,t+\vv}^\mu}, \LL_{t+\vv,T}^{\LL_{t,t+\vv}^\mu})
\big(\e^{\int_t^T V(r,X_{r,T}^{x,\mu},\LL_{r,T}^\mu)d r}- \e^{\int_{t+\vv}^T V(r,X_{r,T}^{x,\mu},\LL_{r,T}^\mu)d r}\big)\Big],\\
&I_2(\vv):= \E \bigg[\int_t^{t+\vv} F(r,X_{t,r}^{x,\mu}, \LL_{t,r}^\mu)\e^{\int_t^r V(\theta, X_{t,\theta}^{x,\mu},\LL_{t,\theta}^\mu)\d\theta}\d r\bigg],\\
&I_3(\vv) :=   \E \bigg[\int_{t+\vv}^T   F(r,X_{t,r}^{x,\mu}, \LL_{t,r}^\mu)\big(\e^{\int_t^r V(\theta, X_{t,\theta}^{x,\mu},\LL_{t,\theta}^\mu)\d\theta}- \e^{\int_{t+\vv}^r V(\theta, X_{t,\theta}^{x,\mu},\LL_{t,\theta}^\mu)\d\theta}\big)\d r\bigg].
\end{align*} Therefore,
\beq\label{LMT}\beg{split} & \lim_{\vv\to 0} \ff{U(t,x,\mu)- \E [U(t+\vv, X_{t,t+\vv}^{x,\mu},\LL_{t,t+\vv}^x) ]}\vv\\
&= V(t,x,\mu)
 \E\big[ \Phi (X_{t,T}^{x,\mu},\LL_{t,T}^\mu)\e^{\int_t^T V(r,X_{t,r}^{x,\mu}, \LL_{t,r}^\mu)\d r}\big]+F(t,x,\mu)\\
 &\qquad + V(t,x,\mu) \E \bigg[\int_{t}^T   F(r,X_{t,r}^{x,\mu}, \LL_{t,r}^\mu) \e^{\int_t^r V(\theta, X_{t,\theta}^{x,\mu},\LL_{t,\theta}^\mu)\d\theta}\bigg]\\
 &= (VU+F)(t,x,\mu). \end{split} \end{equation}
By Proposition \ref{P4.2} below,    $U \in C_b^{0,2,2}([0,T]\times \R^d\times \scr P_2)$ and  $\tt \BL_t U(t,x,\mu)$ is continuous in $(t,x,\mu)$.
Then \eqref{ITO} implies
\beg{align*} \E [U(t+\vv, X_{t,t+\vv}^{x,\mu},\LL_{t,t+\vv}^x) ] =  U(t+\vv,x,\mu) +\E\int_t^{t+\vv} \tt\BL_r  U(r, X_{t,r}^{x,\mu}, \LL_{t,r}^\mu)\d r.\end{align*} Combining this with \eqref{LMT} we arrive at
$$-\pp_t U(t,x,\mu)= \lim_{\vv\to 0}\ff{U(t,x,\mu)-U(t+\vv,x,\mu)}\vv= \tt\BL_t U(t,x,\mu) + (UV+F)(t,x,\mu).$$
Therefore, $U$ solves \eqref{PDE}   with continuous  $\tt \BL_t U.$
\end{proof}

The remainder of this section devotes to the proof of the following result.

\beg{prp}\label{P4.2} Under conditions of Theorem $\ref{T4.1}$ and let $U$ be given in Theorem $\ref{T4.1}(1)$.  Then
$U \in C^{0,2,2}_b([0,T]\times \R^d\times\scr P_2),$ so that $\tt \BL_t U$ is continuous on $[0,T]\times \R^d\times \scr P_2$. \end{prp}

We first introduce some notations which will be used in calculations.
\beg{enumerate} \item[$(a)$]
For   $f \in C^2(\R^d)$,
$$(\nn f(x))v_1:= \<\nn f(x),v_1\>=\nn_{v_1}f(x),\ \ (\nn^2f(x))(v_1,v_2) := \Hess_f (v_1,v_2),\ \ x,v_1,v_2 \in\R^d.$$
\item[$(b)$] For $f\in C^{2}(\scr P_2),$
$$\{D f(\mu)\}\phi  := D_\phi f(\mu) =\int_{\R^d}\<Df(\mu)(x),\phi(x)\>\mu(\d x),\ \ \phi\in L^2(\R^d\to\R^d;\mu).$$
\item[$(c)$] Derivatives of   vector or matrix valued functions are given by those of   component functions. For instance, for $f=(f_{ij})  \in C^{1}(\R^d\times \scr P_2\to\R^l\otimes\R^k)$,
\beg{align*} &  \nn_v f(x,\mu) := \big( \<\nn f_{ij}(x,\mu),v\> \big),\ \  D_\phi  f(x,\mu)  :=\big(D_\phi f_{ij}(x,\mu)\big),\end{align*} where $ x,v\in \R^d,  \mu \in \scr P_2$ and $\phi\in L^2(\R^d\to\R^d;\mu).$
  \end{enumerate}

We will also need the following notion of  uniform  boundedness and continuity.

\beg{defn} Let $\mathbb B$ be a Banach space, and let $E$  be  a topological space. The family
$$\big\{\eta(x)\in L^1(\OO\to \mathbb B;\P):  x\in E\big\}$$ is called  $L^{\infty-}(\P)$    bounded   continuous, if for any $p\ge 1$,
$$\sup_{x\in E}\E\|\eta(x)\|^p<\infty, \ \  \lim_{y\to x} \E\|\eta (x)-\eta(y)\|^p=0,\ \ x\in E.$$\end{defn}

Let  $\L(\mathbb B_1\to\mathbb B_2)$ denote the space of all bounded linear operators from a Banach space $\mathbb B_1$ to the other one $\mathbb B_2$. When $\mathbb B_1$ and $\mathbb B_2$ are finite-dimensional Hilbert spaces,  we regard $\L(\mathbb B_1\to\mathbb B_2)$ as  Euclidean space.
The following   lemma  can be easily proved   by using It\^o's formula, so we omit the proof to save space.

\beg{lem}\label{LP} Let $k,l\ge 1$, and let
\beg{align*} &B_1: \OO\times [0,T]\times \R^l\times\scr P_2\to \R^k,\ \ \Sigma_1:  \OO\times [0,T]\times \R^l\times\scr P_2\to \R^k\otimes \R^m,\\
&B_2: \OO\times [0,T]\times \R^l\times\scr P_2\to \R^k\otimes\R^k,\ \ \Sigma_2:  \OO\times [0,T]\times \R^l\times\scr P_2\to\L(\R^k\to \R^k\otimes \R^m)\end{align*}
 be progressively measurable. If  $\{B_2, \Sigma_2\}$ are uniformly bounded and continuous in $(t, x,\mu)\in [0,T]\times \R^l\times\scr P_2$, and $\{B_1(t,x,\mu), \Sigma_1(t,x,\mu)\}$ are  $L^{\infty-}(\P)$   bounded    continuous, then for any $e\in \R^k$ and $(x,\mu)\in \R^l\times \scr P_2$,
 the solution $(\eta_{s,t}^{x,\mu} )_{t\in [s,T]}$ for the SDE
 $$\d\eta_{s,t}^{x,\mu} = \big\{B_1(t,x,\mu)+ B_2(t,x,\mu)\eta_t^{x,\mu}\big\}\d t +  \big\{\Sigma_1(t,x,\mu)+ \Sigma_2(t,x,\mu)\eta_t^{x,\mu}\big\}\d W_t,\ \eta_{s,s}^{x,\mu}=e, t\in [s,T]$$
is  $L^{\infty-}(\P)$   bounded   continuous.
 \end{lem}

In the following  subsections, we  calculate the first and second  order  derivatives of $(X_{s,t}^{x,\mu},\LL_{s,t}^\mu)$ in $x$ and $\mu$ respectively,  which will be used in  the proof of Proposition \ref{P4.2}.

\subsection{Formulas for $\nn X_{s,t}^{x,\mu}$ and $\nn^2 X_{s,t}^{x,\mu}$}

Let $\{e_i\}_{1\le i\le d}$ be the canonical orthonormal basis of $\R^d$.
Given $(\LL_{s,t}^\mu)_{t\ge s}$,   the SDE \eqref{E1} becomes the classical one  with random coefficients of bounded and continuous first and second order derivatives in $x$.
So,   when $\nn b(t,x,\mu)$ and $\nn \si(t,x,\mu)$ are $L^{\infty-}(\P)$   bounded   continuous,   by taking   $\pp_{x_i}$  to $X_{s,t}^{x,\mu}$ in \eqref{E1}, we see that for any $1\le i\le d$,
$$v_{s,t}^{i,x,\mu}:= \pp_{x_i} X_{s,t}^{x,\mu},\ \ t\ge s $$
solves the linear SDE
 \beq\label{GRT} \beg{split} &\d v_{s,t}^{i,x,\mu} = \Big[\big\{\nn b(t,\cdot, \LL_{s,t}^\mu)(X_{s,t}^{x,\mu})\big\} v_{s,t}^{i,x,\mu} \Big] \d t + \Big[\big\{\nn \si (t,\cdot, \LL_{s,t}^\mu)(X_{s,t}^{x,\mu})\big\} v_{s,t}^{i,x,\mu}\Big]\d W_t,\\
&\ t\ge s,  \ v_{s,s}^{i,x,\mu}= e_i.\end{split} \end{equation}
If moreover $\nn^2 b(t,x,\mu)$ and $\nn^2 \si(t,x,\mu)$ are $L^{\infty-}(\P)$   bounded   continuous, then  by taking $\pp_j$  to the SDE \eqref{GRT},   we see that for   $1\le j\le d$
$$ v_{s,t}^{i,j,x,\mu}:=\pp_{x_i}\pp_{x_j} X_{s,t}^{x,\mu}, \ \ t\ge s$$ solves  the    SDEs
\beg{align*}
&\d v_{s,t}^{i,j,x,\mu}= \Big[ \big\{\nn b(t,\cdot,\LL_{s,t}^\mu)(X_{s,t}^{x,\mu}) \big\} v_{s,t}^{i,j,x,\mu} + \big\{\nn^2 b(t,\cdot,,\LL_{s,t}^\mu)(X_{s,t}^{x,\mu})\big\}(v_{s,t}^{i,x,\mu}, v_{s,t}^{j,x,\mu})\Big]\d t\\
&\quad + \Big[ \big\{\nn \si(t,\cdot,\LL_{s,t}^\mu)(X_{s,t}^{x,\mu}) \big\} v_{s,t}^{i,j,x,\mu} + \big\{\nn^2 \si(t,\cdot,,\LL_{s,t}^\mu)(X_{s,t}^{x,\mu})\big\}(v_{s,t}^{i,x,\mu}, v_{s,t}^{j,x,\mu})\Big]\d W_t,\ \ v_{s,s}^{i,j,x,\mu}=0.\end{align*}
 Combining these with Lemma \ref{LP},  we obtain the following result.

  \beg{lem}\label{LP-1} Assume {\bf (A)} and that $\nn b(t,x,\mu), \nn^2 b(t,x,\mu), \nn \si(t,x,\mu)$ and $\nn^2\si(t,x,\mu)$ are $L^{\infty-}(\P)$   bounded   continuous,
 then so are $\nn X_{s,t}^{x,\mu}$ and $\nn^2 X_{s,t}^{x,\mu}$.        \end{lem}

\subsection{  Formula for $DX_{s,t}^{x,\mu}  $}

We will  establish  the  SDE for $DX_{s,t}^{x,\mu}(y)$ under the following condition {\bf (C)}  on $b$ and $\si$.
\beg{enumerate} \item[{\bf (C)}]  Assume that $b$ and $\si$ are progressively measurable such that the derivatives
  $$ \nn b(t,x,\mu),\ \  \nn \si(t,x,\mu),\ \  D b(t,x,\mu)(y),\ \ D \si(t,x,\mu)(y) $$     are   uniformly  bounded  and continuous
  in $(x,\mu,y)\in \R^d\times \scr P_2\times\R^d.$\end{enumerate}

\beg{lem}\label{LP-2} Assume {\bf (C)}.   Then for any $(x,\mu,y)\in \R^d\times\scr P_2\times\R^d$, $w_{s,t}^{x,\mu}(y):= (D X_{s,t}^{x,\mu})(y)$ for $t\in [s,T]$ exists and solves the  SDE
 \beg{align*}   \d w_{s,t}^{x,\mu}(y) = &\ \bigg[\big\{w_{s,t}^{x,\mu}(y) \big\}^* \nn b(t,\cdot, \LL_{s,t}^\mu)(X_{s,t}^{x,\mu})   +(\nn X_{s,t}^{y,\mu})^*
 \{D b(t, X_{s,t}^{x,\mu}, \cdot)(\LL_{s,t}^\mu)\}(X_{s,t}^{y,\mu}) \\
 & \qquad +  \int_{\R^d} \big\{w_{s,t}^{z,\mu}(y) \big\}^*\big\{Db(t, X_{s,t}^{x,\mu},\cdot)(\LL_{s,t}^\mu)\big\}(X_{s,t}^{z,\mu})\mu(\d z) \bigg] \d t \\
 &+\bigg[\big\{w_{s,t}^{x,\mu}(y) \big\}^* \big\{\nn \si(t,\cdot, \LL_{s,t}^\mu)(X_{s,t}^{x,\mu})\big\}   +(\nn X_{s,t}^{y,\mu})^* \{D \si (t, X_{s,t}^{x,\mu}, \cdot)(\LL_{s,t}^\mu)\}(X_{s,t}^{y,\mu}) \\
 &\qquad + \int_{\R^d} \big\{w_{s,t}^{z,\mu}(y) \big\}^* \big\{D\si (t, X_{s,t}^{x,\mu},\cdot)(\LL_{s,t}^\mu)\big\}(X_{s,t}^{z,\mu})\mu(\d z) \bigg] \d W_t,\ \ w_{s,s}^{x,\mu,y}=0,  \end{align*} where $\big\{w_{s,t}^{x,\mu}(y) \big\}^*$ is the transposition of the matrix $w_{s,t}^{x,\mu}(y).$ Consequently, $(DX_{s,t}^{x,\mu})(y)$ is $L^{\infty-}(\P)$   bounded   continuous. \end{lem}

 To prove the existence of $DX_{s,t}^{x,\mu}$, for fixed $\phi\in L^2(\R^d\to \R^d;\mu)$,   let $\mu_\vv=\mu\circ({\rm Id}+\vv\phi)^{-1}$ and consider
 $$\xi_{s,t}^{x,\vv} := \ff {X_{s,t}^{x,\mu_\vv}- X_{s,t}^{x,\mu}}\vv,\ \ \vv\in (0,1), t\in [s,T].$$  We first establish the SDE for
 $D_\phi X_{s,t}^{x,\mu}:=\lim_{\vv\downarrow 0} \xi_{s,t}^{x,\vv}.$  To this end, we need the following lemma.

 \beg{lem}\label{B1} Assume  {\bf (A)} 
 and let   $\tt\xi_{s,t}^{x,\vv}:= \ff{X_{s,t}^{x+\vv\phi(x),\mu_\vv}- X_{s,t}^{x,\mu_\vv}}\vv.$ Then for any $f\in C^{1,1}(\R^d\times\scr P_2)$ with
 $$K_f:= \sup_{(x,\mu)\in \R^d\times \scr P_2} \big(|\nn f(x,\mu)|^2 + \|Df(x,\mu)\|_{L^2(\mu)}^2\big)<\infty,$$
 the process
 \beg{align*} \Xi_{s,t}^{x,\vv}(f):=&\ \ff{f(X_{s,t}^{x,\mu_\vv},\LL_{s,t}^{\mu_\vv})- f(X_{s,t}^{x,\mu},\LL_{s,t}^{\mu})}\vv -\nn_{\xi_{s,t}^{x,\vv}}f(\cdot,\LL_{s,t}^\mu)(X_{s,t}^{x,\mu}) \\
 &-\int_{\R^d}
 \big\<\xi_{s,t}^{z,\vv}+\tt \xi_{s,t}^{z,\vv}, \{Df(X_{s,t}^{z,\mu},\cdot)(\LL_{s,t}^\mu)\}(X_{s,t}^{z,\mu})\big\>\mu(\d z),\ \ t\in [s,T]\end{align*}
 satisfies
 \beq\label{AD1}  \big|\Xi_{s,t}^{x,\vv}(f)\big|^2\le 8 K_f  \big(  |\xi_{s,t}^{x,\vv}|^2+\mu(|\xi_{s,t}^{\cdot,\vv}+\tt \xi_{s,t}^{\cdot,\vv}|^2)\big),\ \ t\in [s,T], \end{equation}
 \beq\label{AD2} \lim_{\vv\downarrow 0}   \E \big|\Xi_{s,t}^{x,\vv}(f)\big|^2 =0.\end{equation}
 \end{lem}

 \beg{proof}   Let $\eta_r^x= X_{s,t}^{x,\mu}+r(X_{s,t}^{x+\vv\phi(x), \mu_\vv} - X_{s,t}^{x,\mu}),\ r\in [0,1].$ Then $\eta_0^x= X_{s,t}^{x,\mu}, \eta_1^x= X_{s,t}^{x+\vv\phi(x), \mu_\vv}$, so that
 $$\L_{\eta_0|\mu}:= \mu\circ (X_{s,t}^{\cdot,\mu})^{-1}=\LL_{s,t}^\mu,\ \ \L_{\eta_1|\mu}:= \mu\circ(X_{s,t}^{\cdot+\vv\phi,\mu_\vv})^{-1} =\mu_\vv\circ (X_{s,t}^{\cdot,\mu_\vv})^{-1}= \LL_{s,t}^{\mu_\vv}.$$ Moreover,
 $ \ff{\d}{\d r}\eta_r^x= \xi_{s,t}^{x,\vv}+ \tt\xi_{s,t}^{x,\vv}.$
 Then by Lemma \ref{L2.0}, we have
\beg{align*} &\ff{d}{\d r} f(y,\L_{\eta_r|\mu})=\Big\<Df(y,\cdot)(\L_{\eta_r|\mu})(\eta_r), \ff{\d}{\d r}\eta_r\Big\>_{L^2(\mu)} \\
&=\vv \int_{\R^d} \big\<D f(y,\cdot)(\L_{\eta_r|\mu})(\eta_r^z), \xi_{s,t}^{z,\vv}+ \tt\xi_{s,t}^{z,\vv}\big\> \mu(\d z),\ \ \ r\in [0,1], y\in\R^d.\end{align*}
 So, letting $\zeta_{r}^x= (1-r)X_{s,t}^{x,\mu}+rX_{s,t}^{x,\mu_\vv},$ we obtain
 \beg{align*} &\ff{f(X_{s,t}^{x,\mu_\vv},\LL_{s,t}^{\mu_\vv})- f(X_{s,t}^{x,\mu},\LL_{s,t}^{\mu})}\vv
  = \ff 1 \vv \int_0^1 \Big\{\ff{\d}{\d r} f(\zeta_r^x, \L_{\eta_r|\mu})\Big\}\d r \\
 &= \int_0^1 \bigg\{\big\<\nn f(\cdot, \L_{\eta_r|\mu})(\zeta_r^x), \xi_{s,t}^{x,\vv}\big\> + \int_{\R^d}  \big\<D f(\zeta_r^x,\cdot)(\L_{\eta_r|\mu})(\eta_r^z), \xi_{s,t}^{z,\vv}+\tt \xi_{s,t}^{z,\vv}\big\> \mu(\d z)\bigg\}\d r.\end{align*}
 This together with   the definition of $\Xi_{s,t}^{x,\vv}(f)$ gives
 \beq\label{P0} \beg{split}  & \big|\Xi_{s,t}^{x,\vv}(f)\big|^2 = \ \bigg|\int_0^1 \bigg\{\big\<\nn f(\cdot, \L_{\eta_r|\mu})(\zeta_r^x)-\nn f(\cdot,\LL_{s,t}^\mu)(X_{s,t}^{x,\mu}), \xi_{s,t}^{x,\vv}\big\> \\ &  +  \int_{\R^d} \big\<D f(\eta_r^x,\cdot)(\L_{\eta_r|\mu})(\zeta_r^z)- Df(X_{s,t}^{x,\mu},\cdot)(\LL_{s,t}^{\mu})(X_{s,t}^{z,\mu}), \xi_{s,t}^{z,\vv}+\tt \xi_{s,t}^{z,\vv}\big\>\mu(\d z)  \bigg\}\d r \bigg|^2\\
  \le & 8K_f \big(|\xi_{s,t}^{x,\vv}|^2+\mu(|\xi_{s,t}^{\cdot,\vv}+ \tt \xi_{s,t}^{\cdot,\vv}|^2) \big),\end{split} \end{equation}
  which implies \eqref{AD1}. On the other hand,  it is easy to see that   \eqref{UPP}   implies
 \beq\label{PG1}  \sup_{x\in \R^d, \vv\in (0,1)} \E \Big[\sup_{t\in [s,T]}\big\{ |\xi_{s,t}^{x,\vv}|^2+\mu(|\tt\xi_{s,t}^{\cdot,\vv}|^2)\big\}\Big]\le c\mu(|\phi|^2),\ \ \phi\in L^2(\R^d\to\R^d;\mu)\end{equation}
 for some constant $c>0$. Combining this with   the facts that  $(\nn f, Df)$ is bounded continuous,
 $\lim_{r\to 0}\zeta_r^z= X_{s,t}^{z,\mu},$ and $\lim_{r\to 0} \L_{\eta_r|\mu}=\LL_{s,t}^\mu$,  we may apply  the dominated convergence theorem to deduce \eqref{AD2} from the first equality in \eqref{P0} with $\vv\downarrow 0$.
 \end{proof}

 \beg{lem}\label{B2} Assume {\bf (C)}. For any $(s,x,\mu)\in [0,T]\times\R^d,\times\scr P_2$  and   $\phi\in L^2(\R^d\to\R^d)$, $w_{s,t}^{x,\mu,\phi}:= D_\phi X_{s,t}^{x,\mu}$ for $t\in [s,T]$ exists in
 $L^2(\OO\to C([s,T]\to\R^d);\P)$, and there exists a constant $C>0$  such that
 \beq\label{NN1} \E\Big[\sup_{ s\le t\le T} |w_{s,t}^{x,\mu,\phi}|^2\Big]\le C \mu(|\phi|^2),\ \ (s,x,\mu)\in [0,T]\times\R^d\times\scr P_2.\end{equation}    Moreover, for any
  $ t\in [s,T],$
 \beq\label{NN1'} \beg{split}    w_{s,t}^{x,\mu,\phi}
  =   &\ \int_s^t   \Big\{\nn_{w_{s,r}^{x,\mu,\phi} }   b(r, \cdot, \LL_{s,r}^\mu)(X_{s,r}^{x,\mu}) \Big\}\d r + \int_s^t   \Big\{\nn_{w_{s,t}^{x,\mu,\phi}} \si(t, \cdot, \LL_{s,t}^\mu)(X_{s,t}^{x,\mu} )\Big\}\d W_r\\
    &   +   \int_s^t\bigg( \int_{\R^d} \big\<\{D b (r, X_{s,r}^{x,\mu},\cdot)(\LL_{s,r}^\mu)\}(X_{s,r}^{z,\mu}), w_{s,r}^{z,\mu,\phi}+\nn_{\phi(z)}X_{s,t}^{z,\mu}\big\>\mu(\d z)  \bigg) \d r \\
  & +\int_s^t \bigg( \int_{\R^d} \big\<\{D \si (r, X_{s,r}^{x,\mu},\cdot)(\LL_{s,t}^\mu)\}(X_{s,r}^{z,\mu}), w_{s,r}^{z,\mu,\phi}+ \nn_{\phi(z)}X_{s,t}^{z,\mu}\big\>\mu(\d z)\bigg)\d W_r.\end{split}\end{equation}
 \end{lem}

 \beg{proof} To prove the existence of $w_{s,t}^{x,\mu,\phi}:= D_\phi X_{s,t}^{x,\mu}$   in
 $L^2(\OO\to C([s,T]\to\R^d);\P)$, it suffices to show
 \beq\label{N1} \lim_{\vv,\dd\downarrow 0} \E\Big[\sup_{t\in [s,T]} |\xi_{s,t}^{x,\vv}-\xi_{s,t}^{x,\dd}|^2\Big]=0.\end{equation}
 By the definition of $\xi_{s,t}^{x,\vv}$ and letting
 $$\Xi_{s,t}^{x,\vv}(b)= \big(\Xi_{s,t}^{x,\vv}(b_i)\big)_{1\le i\le d}, \ \ \Xi_{s,t}^{x,\vv}(\si)= \big(\Xi_{s,t}^{x,\vv}(\si_{i,j})\big)_{1\le i\le d,1\le j\le m},$$
 we obtain
 \beq\label{N0} \beg{split} &\xi_{s,t}^{x,\vv}  = \ff 1 \vv \int_s^t \big\{b(r,X_{s,r}^{x,\mu_\vv}, \LL_{s,r}^{\mu_\vv}) - b(r,X_{s,r}^{x,\mu}, \LL_{s,r}^{\mu})\big\}\d r \\
&\qquad \qquad  \qquad +\ff 1 \vv\int_s^t  \big\{\si(r,X_{s,r}^{x,\mu_\vv}, \LL_{s,r}^{\mu_\vv}) - \si(r,X_{s,r}^{x,\mu}, \LL_{s,r}^{\mu})\big\}\d W_r  \\
 &= \int_s^t \Big\{\Xi_{s,r}^{x,\vv}(b) +\nn_{\xi_{s,r}^{x,\vv}} b(r,\cdot,\LL_{s,r}^\mu)(X_{s,r}^{x,\mu})\Big\}  \d r \\
  &+\int_s^t \bigg\{  \int_{\R^d} \big\<\{Db(r,X_{s,r}^{x,\mu}, \cdot)(\LL_{s,t}^\mu)\}(X_{s,t}^{z,\mu}), \xi_{s,r}^{z,\vv}+\tt  \xi_{s,r}^{z,\vv}\big\>\mu(\d z) \bigg\}\d r \\
& +\int_s^t\Big\{ \Xi_{s,r}^{x,\vv}(\si) + \nn_{\xi_{s,r}^{x,\vv}} \si(r,\cdot,\LL_{s,r}^\mu)(X_{s,r}^{x,\mu}) \Big\}\d W_r\\
&+ \int_s^t \bigg\{  \int_{\R^d} \big\<\{D\si(r,X_{s,r}^{x,\mu}, \cdot)(\LL_{s,t}^\mu)\}(X_{s,t}^{z,\mu}), \xi_{s,r}^{z,\vv}+\tt  \xi_{s,r}^{z,\vv}\big\>\mu(\d z) \bigg\}\d W_r.\end{split} \end{equation}  Combining this with {\bf (C)} and using the BDG inequality, we may find out a constant $C>0$ such that
 for any $t\in [s,T]$,
 \beq\label{N2}\beg{split}   \E\Big[\sup_{r\in [s,t]} |\xi_{s,r}^{x,\vv}-\xi_{s,r}^{x,\dd}|^2\Big]
  \le &\ C \E\int_s^t \Big\{|\Xi_{s,r}^{x,\vv}(b)- \Xi_{s,r}^{x,\dd}(b)|^2 +\|\Xi_{s,r}^{x,\vv}(\si)- \Xi_{s,r}^{x,\dd}(\si)\|^2\\
 &\qquad  + |\xi_{s,r}^{x,\vv}-\xi_{s,r}^{x,\dd}|^2 + \mu(|\xi_{s,r}^{\cdot,\vv}-\xi_{s,r}^{\cdot,\dd}|^2+|\tt\xi_{s,r}^{\cdot,\vv}-\tt\xi_{s,r}^{\cdot,\dd}|^2)\Big\}\d r.\end{split} \end{equation}
 Integrating both sides with respect to $\mu(\d x)$, we obtain
 \beg{align*} \E\mu(|\xi_{s,t}^{\cdot,\vv}-\xi_{s,t}^{\cdot,\dd}|^2)\le &\ C   \E\int_s^t \mu\big(|\Xi_{s,r}^{\cdot,\vv}(b)- \Xi_{s,r}^{\cdot,\dd}(b)|^2 +\|\Xi_{s,r}^{\cdot,\vv}(\si)- \Xi_{s,r}^{\cdot,\dd}(\si)\|^2+|\tt\xi_{s,r}^{\cdot,\vv}-\tt\xi_{s,r}^{\cdot,\dd}|^2\big)\d r\\
 &+ 2 C\int_s^t \E   \mu(|\xi_{s,r}^{\cdot,\vv}-\xi_{s,r}^{\cdot,\dd}|^2) \d r,\ \ t\in [s,T].\end{align*} Then by Grownwall's inequality,
  \eqref{AD2}, \eqref{PG1}, and the existence   of
 $$\lim_{\vv\downarrow 0} \tt\xi_{s,r}^{\cdot,\vv} =\nn_{\phi} X_{s,t}^{\cdot,\mu}\   \text{in}\ L^2(\P)$$    as explained in Subsection 4.1, which implies
 $\lim_{\vv,\dd\downarrow 0} \E|\tt\xi_{s,r}^{\cdot,\vv}-\tt\xi_{s,r}^{\cdot,\dd}|^2=0,$     we derive
\beg{align*} & \lim_{\vv,\dd\downarrow 0} \sup_{t\in [s,T]} \E \mu(|\xi_{s,t}^{\cdot,\vv}-\xi_{s,t}^{\cdot,\dd}|^2)\\
&\le C \e^{2CT} \lim_{\vv,\dd\downarrow 0}  \E\int_s^T \mu\big(|\Xi_{s,r}^{\cdot,\vv}(b)- \Xi_{s,r}^{\cdot,\dd}(b)|^2 +\|\Xi_{s,r}^{\cdot,\vv}(\si)- \Xi_{s,r}^{\cdot,\dd}(\si)\|^2
+|\tt\xi_{s,r}^{\cdot,\vv}-\tt\xi_{s,r}^{\cdot,\dd}|^2\big)\d r=0.\end{align*} Substituting this into \eqref{N2} and using Gronwall's inequality again,  we arrive at
\beg{align*}   &\lim_{\vv,\dd\downarrow 0} \E \Big[\sup_{t\in [s,T]} |\xi_{s,t}^{x,\vv}-\xi_{s,t}^{x,\dd}|^2\Big]\\
 &\le C\e^{CT} \lim_{\vv,\dd\downarrow 0}   \E\int_s^T \Big\{|\Xi_{s,r}^{x,\vv}(b)- \Xi_{s,r}^{x,\dd}(b)|^2 +\|\Xi_{s,r}^{x,\vv}(\si)- \Xi_{s,r}^{x,\dd}(\si)\|^2  \\
 &\qquad\qquad \qquad\qquad\qquad\qquad\qquad  +\mu\big(|\xi_{s,r}^{\cdot,\vv}- \xi_{s,r}^{\cdot,\dd}|^2+|\tt\xi_{s,r}^{\cdot,\vv}-\tt\xi_{s,r}^{\cdot,\dd}|^2\big)\Big\}\d r=0.\end{align*} Therefore, \eqref{N1} holds, so that $$w_{s,t}^{x,\mu,\phi}:= D_\phi X_{s,t}^{x,\mu}=\lim_{\vv\downarrow 0} \xi_{s,t}^{x,\vv},\ \ \ t\in [s,T]$$   exists in
 $L^2(\OO\to C([s,T]\to\R^d);\P),$   and \eqref{NN1} follows from \eqref{PG1}. Moreover, by   {\bf (C)} and Lemma \ref{B1}, we may let $\vv\downarrow 0$ in \eqref{N0}
to derive the desired equation for $w_{s,t}^{x,\mu,\phi}.$
 \end{proof}

\beg{proof}[Proof of Lemma \ref{LP-2}]   By \eqref{NN1}, $(DX_{s,t}^{x,\mu})_{t\in [s,T]}$ exists with
\beq\label{*W1} \<DX_{s,t}^{x,\mu},\phi\>_{L^2(\mu)}= D_\phi X_{s,t}^{x,\mu}= w_{s,t}^{x,\mu,\phi},\ \ \phi\in L^2(\R^d\to\R^d;\mu).\end{equation}
On the other hand, let $w_{s,t}^{x,\mu}(y)$ solve the SDE in Lemma \ref{LP-2}. Then $\tt w_{s,t}^{x,\mu,\phi}:= \<w_{s,t}^{x,\mu},\phi\>_{L^2(\mu)}$ solves the SDE in Lemma \ref{B2} for $w_{s,t}^{x,\mu,\phi}$. By the uniqueness, we have $w_{s,t}^{x,\mu,\phi}=\tt w_{s,t}^{x,\mu,\phi}.$ Combining this with \eqref{*W1}, we obtain
$\mu$-a.e. $w_{s,t}^{x,\mu}= DX_{s,t}^{x,\mu}$. Then the proof is finished. \end{proof}

\subsection{ Some other derivatives }
We first present a formula for    $Df(\LL_{s,t}^\mu)$.
\beg{lem}\label{LPP} Assume {\bf (C)}.  For any $f\in C_b^{1}(\scr P_2),$
 \beq\label{YP} \beg{split} & \{Df(\LL_{s,t}^\cdot)(\mu)\}(y)\\
 & =   \big( \nn X_{s,t}^{y,\mu}\big)^* \big\{(Df)(\LL_{s,t}^\mu)\big\}(X_{s,t}^{y,\mu})
  +  \int_{\R^d} \big(DX_{s,t}^{x,\mu}\big)^* (y)\{(Df)(\LL_{s,t}^\mu)\}(X_{s,t}^{x,\mu})\mu(\d x).\end{split}\end{equation} \end{lem}

 \beg{proof} Let $\phi\in L^2(\R^d\to\R^d;\mu)$. Since $\LL_{s,t}^\mu= \mu\circ (X_{s,t}^{\cdot,\mu})^{-1}$, for any $\vv>0$  we have
 \beg{align*} &\int_{\R^d} h(z) \big(\LL_{s,t}^{\mu\circ ({\rm Id}+\vv\phi)^{-1}}\big)(\d z) = \int_{\R^d} h\big(X_{s,t}^{x,\mu\circ ({\rm Id}+\vv\phi)^{-1}}\big)\big(\mu\circ ({\rm Id}+\vv\phi)^{-1}\big)(\d x)\\
 &= \int_{\R^d} h\big(X_{s,t}^{x+\vv\phi(x), \mu\circ ({\rm Id}+\vv\phi)^{-1}}\big)\mu(\d x),\ \ h\in \B_b(\R^d).\end{align*}
 So, $\LL_{s,t}^{\mu\circ ({\rm Id}+\vv\phi)^{-1}}$ is  the law of
 $$x\mapsto X_{s,t}^{x+\vv\phi(x), \mu\circ ({\rm Id}+\vv\phi)^{-1}}$$   on the probability space $(\R^d,\B(\R^d),\mu)$. Therefore, by Lemmas \ref{L2.0} and \ref{LP-2}, we obtain
 \beg{align*} &\<Df(\LL_{s,t}^\cdot)(\mu), \phi\>_{L^2(\mu)} := \ff{\d}{\d\vv} f(\LL_{s,t}^{\mu\circ ({\rm Id}+\vv\phi)^{-1}})\Big|_{\vv=0}\\
 &=  \int_{\R^d} \Big\<\{(Df)(\LL_{s,t}^{\mu}) \} (X_{s,t}^{x,\mu}),\ff{\d}{\d\vv} X_{s,t}^{x+\vv\phi(x),\mu\circ ({\rm Id}+\vv\phi)^{-1}}\Big|_{\vv=0}\Big\> \mu(\d x)\\
 &=\int_{\R^d} \big\< \{(Df)(\LL_{s,t}^{\mu}) \} (X_{s,t}^{x,\mu}), \nn_{\phi(x)}X_{s,t}^{x,\mu} +  D_\phi X_{s,t}^{x,\mu}\big\>\mu(\d x)\\
 &= \int_{\R^d} \big\<(\nn X_{s,t}^{x,\mu})^* \{(Df)(\LL_{s,t}^\mu) \}(X_{s,t}^{x,\mu}), \phi(x)\big\>\mu(\d x) \\
 &\qquad+ \int_{\R^d\times\R^d} \big\<(D X_{s,t}^{x,\mu})^*(y) \{(Df)(\LL_{s,t}^{x,\mu})\}(X_{s,t}^{x,\mu}), \phi(y)\big\>\mu(\d x)\mu(\d y)\\
 &= \bigg\< \big( \nn X_{s,t}^{\cdot,\mu}\big)^* \big\{(Df)(\LL_{s,t}^\mu)\big\}(X_{s,t}^{\cdot,\mu})
  +  \int_{\R^d} \big(DX_{s,t}^{x,\mu}\big)^* (\cdot)\{(Df)(\LL_{s,t}^\mu)\}(X_{s,t}^{x,\mu})\mu(\d x), \phi\bigg\>_{L^2(\mu)}.\end{align*}
   Therefore, \eqref{YP} holds. \end{proof}

Next, when $b,\si\in C^{0,2,2}_b([0,T]\times\R^d\times \scr P_2)$,    by making derivatives to the SDE for $w_{s,t}^{x,\mu}(y)$ presented in  Lemma \ref{LP-2},
  we derive the following result.

\beg{lem}\label{LP-3} Assume that $b,\si\in C^{0,2,2}_b([0,T]\times\R^d\times \scr P_2).$
Then  all derivatives $$  \{D\nn X_{s,t}^{x,\mu}\}(y),\ \nn \{D X_{s,t}^{\cdot,\mu}(y)\}(x),
\ \nn \{DX_{s,t}^{y,\mu}(\cdot)\}(y),\ D^2X_{s,t}^{x,\mu}(y,z)$$ are $L^{\infty-}(\P)$   bounded   continuous.  \end{lem}

\beg{proof}  (a) We first consider $\{D\nn X_{s,t}^{x,\mu}\}(y)$. Since $b,\si\in C^{0,2,2}_b([0,T]\times\R^d\times \scr P_2)$, by \eqref{GRT} and Lemmas \ref{LP-1}-\ref{LP-2},
$v_{s,t}^{x,\mu}:=\nn_v X_{s,t}^{x,\mu}$ for $v\in \R^d$ solves the SDE
$$\d v_{s,t}^{x,\mu} = Z_1(t,x,\mu) v_{s,t}^{x,\mu}\d t+ \{Z_2(t,x,\mu) v_{s,t}^{x,\mu}\}\d W_t,\ \ v_{s,s}^{x,\mu}= v,$$
where $$Z_1:  [0,T]\times\R^d\times \scr P_2\to \R^d,\ \ Z_2: [0,T]\times\R^d\times \scr P_2\to \R^d\otimes\R^d$$   are progressively measurable and satisfy
\beg{enumerate} \item[{\bf (D)}] $Z_1(t,x,\mu)$ and $Z_2(t,x,\mu)$ are uniformly bounded and continuous in $(t,x,\mu)\in [0,T]\times \R^d\times\scr P_2$;    $D Z_1(t,x,\mu)(y)$ and $D Z_2(t,x,\mu)(y)$
are $L^{\infty-}(\P)$  bounded   continuous.\end{enumerate}
  Then for any $\phi\in L^2(\R^d\to\R^d;\mu)$ and $\mu_\vv:= \mu\circ({\rm Id} +\vv\phi)^{-1}$ for small $\vv>0$, $\gg_{s,t}^\vv:= \ff{v_{s,t}^{x,\mu_\vv} - v_{s,t}^{x,\mu}}\vv$ solves the SDE
\beg{align*} &\d \gg_{s,t}^\vv=   \{Z_1(t,x,\mu) \gg_{s,t}^{\vv}\}\d t+ \{Z_2(t,x,\mu) \gg_{s,t}^{\vv}\}\d W_t\\
&+ \ff{\{Z_1(t,x,\mu^\vv) -Z_1(t,x,\mu)\}v_{s,t}^{x,\mu_\vv}}\vv\, \d t+ \ff{\{Z_2(t,x,\mu_\vv) -Z_2(t,x,\mu)\}  v_{s,t}^{x,\mu_\vv}}\vv\,\d W_t,\ \ \eta_{s,s}^\vv=0.\end{align*}
By {\bf (D)}, we may repeat the proof of Lemma \ref{B2} to conclude that  $D_\phi v_{s,t}^{x,\mu} := \lim_{\vv\downarrow 0} \eta_{s,t}^\vv$ exists and solves the SDE
\beg{align*} \d \{D_\phi v_{s,t}^{x,\mu}\}= &\ \big\{Z_1(t,x,\mu) D_\phi v_{s,t}^{x,\mu} + (D_\phi Z_1(t,x,\mu))v_{s,t}^{x,\mu}\big\}\d t\\
& + \big\{Z_2(t,x,\mu)D_\phi v_{s,t}^{x,\mu} + (D_\phi Z_2(t,x,\mu))v_{s,t}^{x,\mu} \big\}  \d W_t,\ \ D_\phi v_{s,s}^{x,\mu}=0.\end{align*}
Hence,  $D v_{s,t}^{x,\mu}(y)$ solves the SDE
\beg{align*} \d \{Dv_{s,t}^{x,\mu}(y)\}= &\ \big\{Z_1(t,x,\mu) D  v_{s,t}^{x,\mu} (y)+ (D  Z_1(t,x,\mu)(y)) v_{s,t}^{x,\mu}\big\}\d t\\
&+ \big\{Z_2(t,x,\mu)D  v_{s,t}^{x,\mu} (y)+ (D  Z_2(t,x,\mu)(y))v_{s,t}^{x,\mu} \big\}  \d W_t,\ \ D v_{s,s}^{x,\mu}(y)=0.\end{align*}
Therefore, by Lemma \ref{LP-1} and {\bf (D)}, Lemma \ref{LP}  yields  that $\{D\nn X_{s,t}^{x,\mu}\}(y)$ is $L^{\infty-}(\P)$  bounded   continuous.

(b) To calculate $ \nn\{DX_{s,t}^{\cdot,\mu}(y)\}(x), \nn \{D X_{s,t}^{x,\mu}(\cdot)\}(y)$ and $D^2X_{s,t}^{x,\mu}(y,z):= D\{DX_{s,t}^{x,\mu}(y)\}(z),$ we reformulate the SDE in Lemma \ref{LP-2} for $w_{s,t}^{x,\mu}(y):= D X_{s,t}^{x,\mu}(y)$ as
$$\d w_{s,t}^{x,\mu} = \big\{A_1 (t,x,\mu) w_{s,t}^{x,\mu} + A_2(t,x,\mu)\big\} \d t+ \{B_1(t,x,\mu) w_{s,t}^{x,\mu}+ B_2(t,x,\mu)\big\}\d W_t,\ \ w_{s,s}^{x,\mu}= 0,$$
where, due to Lemmas \ref{LP-1}-\ref{LP-2} and (a),
$\{A_i,B_i\}_{i=1,2}$ are progressively  measurable maps such that
\beg{enumerate} \item[$\bullet$] $ A_1$ and $ B_1$ are uniformly bounded and continuous in $(t,x,\mu)\in [0,T]\times\R^d\times \scr P_2$;
\item[$\bullet$]   $\{A_i, B_i, \nn A_i, \nn B_i, D A_i, D B_i\}_{i=1,2}$ are  $L^{\infty-}(\P)$  bounded  continuous in corresponding arguments.\end{enumerate}
So, as explained in (a),  by taking derivatives $\pp_{x_i}, \pp_{y_i}$ and $D_\phi$ to this SDE respectively and applying Lemma \ref{LP},   we prove  that $\pp_{y_i} D X_{s,t}^{x,\mu} (y)$ and $D^2X_{s,t}^{x,\mu}(y,z)$
are $L^{\infty-}(\P)$   bounded   continuous in related arguments. We omit the details to save space.
 \end{proof}

\subsection{Proof of  Proposition \ref{P4.2}}

Since $b,\si\in C^{0,2,2}_b([0,T]\times\R^d\times \scr P_2)$,  assertions in Lemmas \ref{LP-1},   \ref{LP-2}, and   \ref{LP-3} hold. Then  it is straightforward to show that $U$ given in Theorem \ref{T4.1}(1) is in the class $C^{0,2,2}([0,T]\times\scr P_2)$.

Firstly, for any $1\le i\le d$, by taking derivative $\pp_{x_i}$ to the formula of $U$, we obtain
\beg{align*} &\pp_{x_i} U(t,x,\mu) =  \E\Big[\big\<\nn\Phi(\cdot, \LL_{t,T}^\mu)(X_{t,T}^{x,\mu}), \pp_{x_i} X_{t,T}^{x,\mu}\big\> \e^{\int_t^T V(r,X_{t,r}^{x,\mu}, \LL_{t,r}^\mu)\d r}\Big]\\
&+\E\bigg[\Phi(X_{t,T}^{x,\mu},\LL_{t,T}^\mu)  \e^{\int_t^T V(r,X_{t,r}^{x,\mu}, \LL_{t,r}^\mu)\d r} \int_t^T\big\<\nn V(r,\cdot, \LL_{t,r}^\mu)(X_{t,r}^{x,\mu}), \pp_{x_i} X_{t,r}^{x,\mu}\big\>\d r\bigg]\\
&+\E\int_t^T  \big\<\nn F(r,\cdot, \LL_{t,r}^\mu)(X_{t,r}^{x,\mu}), \pp_{x_i} X_{t,r}^{x,\mu}\big\> \e^{\int_t^r V(\theta,X_{t,\theta}^{x,\mu}, \LL_{t,\theta}^\mu)\d \theta}\,\d r \\
&+\E\int_t^T\bigg\{F(r,X_{t,r}^{x,\mu},\LL_{t,r}^\mu)  \e^{\int_t^r V(\theta,X_{t,\theta}^{x,\mu}, \LL_{t,\theta}^\mu)\d \theta} \int_t^r\big\<\nn V(\theta,\cdot, \LL_{t,\theta}^\mu)(X_{t,\theta}^{x,\mu}), \pp_{x_i} X_{t,\theta}^{x,\mu}\big\>\d \theta\bigg\}\,\d r.\end{align*}
By assumptions on $\Phi, V, F$ and  Lemmas \ref{LP-1}, \ref{LP-2}  and \ref{LP-3},   this formula implies that $\nn U(t,x,\mu)$ is bounded and continuous. Moreover,
  by taking derivatives $\pp_{x_j}$ and $D$ to the formula, we conclude that $ \nn^2 U(t,x,\mu)$
 and  $D \{\nn X_{s,t}^{x,\mu}\}(y)$ are bounded and continuous as well.

 Similarly, we may prove the assertion for $DU(t,x,\mu)(y), \pp_{x_i} \{D U(t,x,\mu)(y)\}, \pp_{y_i}  \{D U(t,x,\mu)(y)\}$ and $D^2 U(t,x,\mu)(y,z)$. For simplicity, we only consider the case for $V=F=0$,
 for the general case the formulation is only more complicated due to derivatives to $F$ and $V$, but there is no any essential difference for the proof.
For $V=F=0$ the formula for $U$ becomes
$$U(t,x,\mu)= \E \Phi(X_{t,T}^{x,\mu},\LL_{s,t}^\mu).$$
 Then by  \eqref{YP} and the   chain rule  we obtain
\beg{align*}  D U(t,x,\mu)(y)
= &\ \E\bigg[ \big\{\nn\Phi(\cdot,\LL_{s,t}^{\mu})(X_{s,t}^{x,\mu}) \big\}(DX_{s,t}^{x,\mu})(y)
    +(\nn X_{s,t}^{y,\mu})^* \big\{(D\Phi(t,X_{s,t}^{x,\mu}, \cdot)(\LL_{s,t}^\mu)\big\}(X_{s,t}^{y,\mu})\\
 &\qquad   +\int_{\R^d}(DX_{s,t}^{z,\mu})^* (y) \big\{D\Phi(X_{s,t}^{x,\mu},\cdot)(\LL_{s,t}^\mu)\big\}(X_{s,t}^{z,\mu})\mu(\d z) \bigg].
 \end{align*}
 Since $\Phi\in C_b^{2,2}(\R^d\times\scr P_2)$,  by Lemmas \ref{LP-1}, \ref{LP-2}  and \ref{LP-3} we deduce from this formula that $D U(t,x,\mu)(y)$ is bounded and continuous. Moreover,
 by taking derivatives $  \pp_{x_i},\pp_{y_i}, D $ to this formula, we conclude that $\pp_{x_i} \{D U(t,x,\mu)(y)\}, \pp_{y_i}  \{D U(t,x,\mu)(y)\}$ and $D^2 U(t,x,\mu)(y,z)$
are bounded and continuous as well. In conclusion, $U\in C^{0,2,2}([0,T]\times\R^d\times\scr P_2).$

\section{Ergodicity  and structure of invariant measures}

In this part, we assume that $b(t,x,\mu)= b(x,\mu)$ and $\si(t,x,\mu)= \si(x,\mu)$ are deterministic, and consider the ergodicity of the diffusion processes generated by $\BL$ and $\tt \BL$.

Recall that a Markov process is called ergodic, if   for any initial distribution, when $t\to\infty$  the process converges weakly to the unique invariant probability measure. For square   integrable Markov processes, the weak convergence is equivalent to the convergence under the Wasserstein distance.
To estimate the Wasserstein distance for solutions to the image SDE \eqref{E1},
we take  the following hypothesis:
\beg{enumerate}\item[{\bf (H)}] $b(t,x,\mu)= b(x,\mu)$ and $\si(t,x,\mu)= \si(x,\mu)$ are deterministic, continuous in $(x,\mu)$ and  do not depend on $t$. There exist constants $\ll\in\R$ and $\kk,\dd,K\ge 0$   such that
\beg{align*} &2 \<b(x,\mu)-b(y,\nu), x-y\>+ \|\si(x,\mu)- \si(y,\nu)\|_{HS}^2\le \kk\W(\mu,\nu)^2-\ll |x-y|^2,\\
&\|\si(x,\mu)- \si(y,\nu)\|_{HS}^2\le K \big\{\W(\mu,\nu)^2+ |x-y|^2\big\},\\
& |b(x,\mu)|^2+\|\si(x,\mu)\|_{HS}^2\le \dd(1+|x|^2+\|\mu\|_2^2),\ \ x,y\in \R^d, \mu,\nu\in \scr P_2.\end{align*}
\end{enumerate}

By Theorem \ref{T2.1},   {\bf (H)} implies the well-posedness of  \eqref{E1}. In the present time-homogenous case, we only consider the solution from time $s=0$,
i.e. $(X_t^{x,\mu}, \LL_t^\mu):= (X_{0,t}^{x,\mu}, \LL_{0,t}^\mu)$ for $t\ge 0$.

Let $P_t(\mu;\cdot)$ and $\tt P_t(x,\mu;\cdot)$ denote the laws of
$\LL_t^\mu$ and $(X_t^{x,\mu}, \LL_t^\mu)$ respectively. Then the associated Markov semigroups $P_t$ and $\tt P_t$ are given by
\beg{align*} &P_t f(\mu):= \E f(\LL_t^\mu)=\int_{\scr P_2} f(\nu)P_t(\mu;\d\nu),\ \ f\in \B_b(\scr P_2),\\
&\tt P_t  g(x,\mu) :=\E g(X_{t}^{x,\mu},\LL_{t}^\mu)=\int_{\R^d\times\scr P_2} g(y,\nu) \tt P_t(x,\mu;\d y,\d\nu), \ \  g\in \B_b(\R^d\times\scr P_2).\end{align*}

Let $\scr P_2(\scr P_2)$ (resp. $\scr P_2(\R^d\times\scr P_2)$) be the set of probability measures on $\scr P_2$ (resp. $\R^d\times\scr P_2$) with finite second moments, and let
${\bf W}_2^{\scr P_2}$ be the $L^2$-Warsserstein distance    on $\scr P_2(\scr P_2)$ induced by $\W_2$, while ${\bf W}_2^{\R^d\times\scr P_2}$ be that
on $\scr P_2(\R^d\times\scr P_2)$ induced by the metric
$$\rr((x,\mu),(y,\nu)):= \ss{|x-y|^2+\W_2(\mu,\nu)^2}.$$

For any $Q\in \scr P_2(\scr P_2)$ and $\tt Q\in \scr P_2(\R^d\times\scr P_2)$, let
$$QP_t =\int_{\scr P_2} P_t (\mu;\cdot)  Q(\d\mu),\ \  \tt Q \tt P_t =\int_{\R^d\times \scr P_2} \tt P_t(x,\mu;\cdot) \tt Q(\d x,\d\mu).$$

In the following two subsections, we first investigate the exponential ergodicity  of the diffusion processes generated by $\BL$ and $\tt\BL$, then figure out the structure of  the invariant probability measures.

\subsection{Exponential ergodicity}

\beg{thm}\label{T3.1} Assume {\bf (H)}.
Then for any $(x,\mu)\in \R^d\times\scr P_2,$
\beq\label{LLPT0} \E \W_2(\LL_t^\mu,\LL_t^\nu)^2\le \W_2(\mu,\nu)^2\e^{-(\ll-\kk)t},\ \ t\ge 0,\end{equation}
\beq\label{LLPT1} \E |X_t^{x,\mu}-X_t^{y,\nu}|^2\le  |x-y|^2\e^{-\ll t} + \W_2(\mu,\nu)^2\e^{-(\ll-\kk)t},\ \ t\ge 0.\end{equation}
Consequently, if $\ll>\kk$ then:
\beg{enumerate}
\item[$(1)$]     $\tt P_t $ has a unique invariant probability measure $\tt \Pi\in \scr P_2(\R^d\times\scr P_2)$  such that for any $\tt Q\in \scr P_2(\R^d\times\scr P_2)$,
\beq\label{ES1}  {\bf W}_2^{\R^d\times\scr P_2}(\tt Q \tt P_t,  \tt \Pi)^2\le 2\e^{-(\ll-\kk)t}{\bf W}_2^{\R^d\times\scr P_2}(\tt Q,\tt\Pi)^2,\ \ t\ge0;\end{equation}
\item[$(2)$] $\Pi:= \tt\Pi(\R^d\times\cdot)$ is the unique invariant probability measure of $P_t$ such that for any $Q\in \scr P_2(\scr P_2)$,
\beq\label{ES2}  {\bf W}_2^{\scr P_2}(Q P_t(\mu;\cdot),  \Pi)^2\le \e^{-(\ll-\kk)t}   {\bf W}_2^{\scr P_2}(Q,  \Pi)^2,\ \ t\ge0.\end{equation}
\end{enumerate}  \end{thm}

\beg{proof} (a) We first prove \eqref{LLPT0} and \eqref{LLPT1}.
Let $\pi\in \scr C(\mu,\nu)$ such that
$$\W_2(\mu,\nu)^2=\int_{\R^d\times\R^d} |x-y|^2\pi(\d x,\d y).$$ Then for any $t\ge 0$,
$$\pi_t:=\pi\circ(X_t^{\cdot,\mu}, X_t^{\cdot,\nu})^{-1}\in \scr C(\LL_t^\mu,\LL_t^\nu), $$
so that
  \beq\label{DR1} \W_2(\LL_t^\mu,\LL_t^\nu)^2\le \int_{\R^d\times\R^d} |x-y|^2\pi_t(\d x,\d y)
  =  \int_{\R^d\times\R^d} |X_t^{x,\mu}-X_t^{y,\nu}|^2\pi(\d x,\d y)=:\ell_t.\end{equation}
Combining this with  {\bf (H)} and It\^o's formula, we obtain
$$ \d |X_t^{x,\mu}-X_t^{y,\nu}|^2 \le \big\{\kk \ell_t -\ll |X_t^{x,\mu}-X_t^{y,\nu}|^2\big\}\d t+\d M_t$$ for some martingale $M_t$, which implies
\beq\label{DR2} \e^{\ll t} \E |X_t^{x,\mu}-X_t^{y,\nu}|^2\le |x-y|^2 + \kk \int_0^t \e^{\ll s} \E \ell_s\d s,\ \ t\ge 0.\end{equation}  Integrating with respect to $\pi(\d x,\d y)$ gives
$$\e^{\ll t} \E\ell_t\le \W_2(\mu,\nu)^2 + \kk \int_0^t \e^{\ell s}  \E\ell_s\d s,\ \ t\ge 0,$$ which together with   Grownwall's lemma and \eqref{DR1} leads to
$$\E\W_2(\LL_t^\mu,\LL_t^\nu)^2\le \E \ell_t\le \W_2(\mu,\nu)^2  \e^{-(\ll -\kk)t},\ \ t\ge 0.$$ Thus, \eqref{LLPT0} holds.
Substituting \eqref{LLPT0} into \eqref{DR2} we arrive at
\beg{align*}  \E |X_t^{x,\mu}-X_t^{y,\nu}|^2&\le \e^{-\ll t}|x-y|^2 + \kk \W_2(\mu,\nu)^2 \e^{-\ll t} \int_0^t \e^{\kk s}\d s\\
&\le \e^{-\ll t}|x-y|^2
+ \W_2(\mu,\nu)^2  \e^{-(\ll -\kk)t}.\end{align*}  Hence, \eqref{LLPT1} holds.

(b) Existence of invariant probability measures.  Consider, for instance $(X_t^{0,\dd_0}, \LL_t^{\dd_0})$, where $\dd_0$ is the Dirac measure at $0\in\R^d.$ Let $\tt\Pi_t=  \tt P_t(0,\dd_0;\cdot)$ be the law of $(X_t^{0,\dd_0}, \LL_t^{\dd_0})$. By the completeness of the Wasserstein space,  if
\beq\label{INV} \lim_{s,t\to\infty} {\bf W}_2^{\R^d\times\scr P_2}(\tt\Pi_t, \tt\Pi_s)^2=0,\end{equation}  then there exists a probability measure $\tt\Pi$ on $\R^d\times \scr P_2$ with
$\|\tt\Pi\|_2^2:= \tt\Pi(\rr^2) <\infty$ such that $\lim_{t\to\infty} {\bf W}_2^{\R^d\times\scr P_2} (\tt\Pi_t,\tt\Pi)=0.$ Consequently, $\tt\Pi$ is an invariant probability measure for $\tt P_t $. Moreover, since the law of $ \LL_t^{\dd_0}$ is $\Pi_t(\R^d\times\cdot),$ which converges to $\Pi:=\tt\Pi(\R^d\times\cdot)$ weakly as $t\to\infty$,  we see that $\Pi$ is an invariant probability measure of $P_t.$

To prove \eqref{INV}, let $t>s\ge 0$. By the Markov property  we have
$$\tt\Pi_t=    P_t(0,\dd_0;\cdot)= \int_{\R^d\times\scr P_2} P_s(x,\mu;\cdot) \tt\Pi_{t-s}(\d x,\d\mu).$$
Combining this with \eqref{LLPT0} and \eqref{LLPT1} we obtain
\beg{equation*}  \beg{split} & {\bf W}_2^{\R^d\times\scr P_2}\big(\tt\Pi_t, \tt\Pi_s)^2\le \int_{\R^d\times\scr P_2} {\bf W}_2^{\R^d\times\scr P_2}(\tt P_s(x,\mu;\cdot), \tt P_s (0,\dd_0;\cdot)\big)^2\,\tt\Pi_{t-s}(\d y,\d \nu)\\
&\le  \int_{\R^d\times\scr P_2}\big\{\E|X_s^{0,\dd_0}-X_s^{x,\mu}|^2 + \W_2(\LL_s^\mu, \LL_s^{\dd_0})^2\big\} \tt\Pi_{t-s}(\d x,\d \mu)\\
&\le \int_{\R^d\times\scr P_2}\big\{|x|^2\e^{-\ll s}  + 2\W_2(\dd_0, \mu)^2\e^{-(\ll-\kk)s} \big\} \tt\Pi_{t-s}(\d x,\d \mu)\\
&= \e^{-\ll s} \E|X_{t-s}^{0,\dd_0}|^2 + 2 \e^{-(\ll-\kk)s} \E \W_2(\dd_0,\LL_{t-s}^{\dd_0})^2= (\e^{-\ll s} + 2  \e^{-(\ll-\kk)s} ) \E|X_{t-s}^{0,\dd_0}|^2.    \end{split} \end{equation*}
So, to prove \eqref{INV} it remains to show that
\beq\label{DR4} \sup_{t\ge 0} \E |X_t^{0,\dd_0}|^2 <\infty.\end{equation}
By assumption {\bf (H)} with $\ll>\kk$, for any $\ll>\ll'>\kk'>\kk$ there exists a constant $c>0$ such that
$$2\<b(x,\mu), x\> + \|\si(x,\mu)\|_{HS}^2\le c+ \kk' \|\mu\|_2^2 -\ll' |x|^2,\ \ (x,\mu)\in \R^d\times\scr P_2.$$
Combining this with   It\^o's formula, and noting that $\|\LL_t^{\dd_0}\|_2^2 = \dd_0(|X_t^{\cdot,\dd_0}|^2)= |X_t^{0,\dd_0}|^2,$ we obtain
$$\d |X_t^{0,\dd_0}|^2 \le \big\{c+(\kk' -\ll')|X_t^{0,\dd_0}|^2\big\}\d t+\d M_t$$ for some martingale $M_t$. This implies
$$\E |X_t^{0,\dd_0}|^2 \le c\int_0^t \e^{-(\ll'-\kk')s}\d s,\ \ t\ge 0.$$ Since $\ll'>\kk'$, we derive    \eqref{DR4} and hence finish the proof of the existence of invariant probability measures. Moreover, the invariant probability measure $\tt\Pi$ satisfies
$$\int_{\R^d\times\scr P_2} (|x|^2+\|\mu\|_2^2)\tt\Pi(\d x,\d\mu)\le \lim_{t\to\infty} \E |X_t^{0,\dd_0}|^2 \le \ff c{\ll'-\kk'}<\infty.$$ Hence, $\tt\Pi\in \scr P_2(\R^d\times\scr P_2).$

(c) It is easy to see that \eqref{ES1} follows from \eqref{LLPT0} and \eqref{LLPT1}.  Indeed,
letting $\GG\in \C(\tt Q,\tt\Pi) $ such that
$${\bf W}_2^{\R^d\times \scr P_2} (\tt Q,\tt\Pi)^2= \int_{(\R^d\times \scr P_2)^2} \rr^2 \d\GG,$$
we deduce from  \eqref{LLPT0}, \eqref{LLPT1}  and $\tt\Pi= \tt \Pi \tt P_t $    that
\beg{align*} &{\bf W}_2^{\R^d\times \scr P_2} (\tt Q \tt P_t, \tt \Pi)^2={\bf W}_2^{\R^d\times \scr P_2} (\tt Q \tt P_t, \tt \Pi \tt P_t)^2 \\
&\le \int_{(\R^d\times\scr P_2)^2} {\bf W}_2^{\R^d\times \scr P_2} (\tt P_t(x,\mu; \cdot), \tt P_t (y,\nu;\cdot))^2\GG(\d x,\d\mu; \d y,\d\nu) \\
&\le   \int_{(\R^d\times\scr P_2)^2}  \E \big\{|X_t^{x,\mu}-X_t^{y,\nu}|^2 +\W_2(\LL_t^\mu, \LL_t^\nu)^2\big\}   \GG(\d x,\d\mu;\d y,\d\nu) \\
&\le \int_{(\R^d\times\scr P_2)^2}    \big\{|x-y|^2\e^{-\ll t}  +2\W_2(\mu,\nu)^2\e^{-(\ll-\kk)t}\big\}  \GG(\d x,\d\mu;\d y,\d\nu)\\
&\le  2\e^{-(\ll-\kk)t} {\bf W}_2^{\R^d\times \scr P_2} (\tt Q, \tt\Pi)^2,\ \ t\ge 0.\end{align*}  In particular, $\tt\Pi$ is the unique invariant probability measure of $P_t$.

(d) As shown  in (b) and (c), \eqref{LLPT0} for $\ll>\kk$ implies that $P_t$ has a unique invariant probability measure $\Pi$ satisfying the estimate \eqref{ES2}. Noting that
$P_t(\mu;\cdot)= \tt P_t(x,\mu; \R^d\times \cdot)$ holds for all $(x,\mu)\in \R^d\times \scr P_2$, we have
$\Pi= \tt\Pi(\R^d\times\cdot).$
\end{proof}

\subsection{Structure of invariant probability measures}

Under condition {\bf (H)}, let $b_0(x)=b(x,\dd_x)$ and $\si_0(x)= \si(x,\dd_x)$. Then the SDE \eqref{E0} is well-posed. Let $P_t^0$ be the associated Markov semigroup.

\beg{thm}\label{T4.2} Assume {\bf (H)}. If $P_t^0$ has an invariant probability measure $\mu_0$, then
$$\tt\Pi_0(\d x,\d\mu):= \mu_0(\d x) \dd_{\dd_x}(\d\mu)$$ is an invariant probability measure of $\tt P_t$. Consequently,
$\Pi_0:= \tt\Pi_0(\R^d\times\cdot)=\int_{\R^d}\dd_{\dd_x}\mu_0(\d x)$ is an invariant probability measure of $P_t$, and
when $\ll>\kk$,
the unique invariant probability measures   $\tt\Pi$ and $\Pi$  in Theorem $\ref{T3.1} $ satisfy  \eqref{E01}. \end{thm}

\beg{proof} Recall that $(X_t^{x,\mu}, \LL_t^\mu) $ solve the SDE
$$\d X_t^{x,\mu}= b(X_t^{x,\mu}, \LL_t^\mu)\d t +\si (X_t^{x,\mu}, \LL_t^\mu)\d W_t,\ \ X_0^{x,\mu}=x,$$
where $\LL_t^\mu:=\mu\circ (X_t^{\cdot,\mu})^{-1}.$ Then, when $\mu=\dd_x$ we have $\LL_t^\mu= \dd_{X_t^{x,\dd_x}}$, so that
$(X_{t}^{x,\dd_x})_{t\ge 0}$ solves the SDE \eqref{E0}.  By the uniqueness of this SDE and that $\mu_0$ is an invariant probability measure of   $P_t^0$,  we obtain
$$\int_{\R^d} \big[\E g(X_t^{x,\dd_x}) \big]\mu_0(\d x) =\int_{\R^d} P_t^0 g(x) \mu_0(\d x)= \int_{\R^d} g(x)\mu_0(\d x),\ \ t\ge 0, g\in \B_b(\R^d).$$
Combining this with $\tt P_t f(x,\dd_x)= \E f(X_t^{x,\dd_x}, \dd_{X_t^{x,\dd_x}})$ for $f\in \B_b(\R^d\times\scr P_2)$, and taking $g(x)=f(x,\dd_x)$, we obtain
\beg{align*} &\int_{\R^d\times\scr P_2} \tt P_t f(x,\mu) \tt\Pi_0(\d x,\d\mu) = \int_{\R^d} \tt P_t f(x,\dd_x) \mu_0(\d x) \\
&=\int_{\R^d} \big[\E f(X_t^{x,\dd_x},\dd_{X_t^{x,\dd_x}})\big]\mu_0(\d x)
 = \int_{\R^d} \big[\E g(X_t^{x,\dd_x})\big]\mu_0(\d x)\\
 &= \int_{\R^d} g(x)\mu_0(\d x)=\int_{\R^d} f(x,\dd_x)\mu_0(\d x)= \int_{\R^d\times\scr P_2} f(x,\mu)\tt\Pi_0(\d x,\d\mu).\end{align*}
Therefore, $\tt\Pi_0$ is an invariant probability measure of $\tt P_t$. In particular, by taking $f(x,\mu)=f(\mu)$, we see that $\Pi_0$ is an invariant probability measure of $P_t$.

Finally, if $\ll>\kk$,  by Theorem \ref{T3.1}, $\Pi$ and $\tt \Pi$ are the unique  invariant probability measures of $P_t$ and $\tt P_t$ respectively.
So,  $\tt \Pi=\tt\Pi_0$ and $\Pi=\Pi_0$; that is, \eqref{E01} holds. \end{proof}

\paragraph{Acknowledgement.} The author is grateful to the referee for valuable suggestions and to Professor Renming Song for helpful conversations.

\end{document}